\documentclass[12pt] {article}
\usepackage{amsmath,amsthm}
\usepackage{amssymb,latexsym}

\usepackage{comment}
\usepackage{tikz-cd}

\title{Non-Archimedean analogue of the space of valuations on convex sets.}
\date{}
\author{Semyon Alesker \footnote{Partially supported by ISF grant 743/22 and  US- Israel BSF grant 2018115.}
\\  { \normalsize Department of Mathematics, Tel Aviv University, Ramat Aviv}
\\  { \normalsize 69978 Tel Aviv, Israel }
\\ {\normalsize e-mail: alesker.semyon75@gmail.com}}

\def\RR{\mathbb{R}}
\def\CC{\mathbb{C}}
\def\QQ{\mathbb{Q}}

\def\ZZ{\mathbb{Z}}

\def\FF{\mathbb{F}}

\def\eps{\varepsilon}
\def\alp{\alpha}
\def\ome{\omega}

\def\lam{\lambda}
\def\Lam{\Lambda}

\def\to{\longrightarrow}

\def\qed { Q.E.D. }

\def\inj{\hookrightarrow}

\usepackage{chngcntr}
\usepackage{enumitem}
\newlist{paragraphlist}{enumerate}{1}
\usepackage{tikz}

\setlist[paragraphlist,1]{leftmargin=*,label={\bfseries \arabic*}}

\counterwithin{paragraphlisti}{subsubsection}

\swapnumbers
\newtheorem{theorem}{Theorem}[section]
\newtheorem{corollary}[theorem]{Corollary}
\newtheorem{lemma}[theorem]{Lemma}
\newtheorem{proposition}[theorem]{Proposition}
\newtheorem{claim}[theorem]{Claim}
\theoremstyle{definition}

\newtheorem{example}[theorem]{Example}
\newtheorem{definition}[theorem]{Definition}
\newtheorem{remark}[theorem]{Remark}
\theoremstyle{conjecture}

\theoremstyle{principle}

  \def\cc{{\cal C}}
\def\cd{{\cal D}}  
  
 \def\ck{{\cal K}} \def\cl{{\cal L}}
\def\cm{{\cal M}}  \def\co{{\cal O}}

 \def\ct{{\cal T}} \def\cu{{\cal U}}
\def\cv{{\cal V}} \def\cw{{\cal W}} \def\cx{{\cal X}}
 
\def\cm{{\cal M}}
\def\pt{\partial}
\def\Vai{Val^\infty}
\numberwithin{equation}{section}
\begin{document}
\maketitle

{\itshape Dedicated to the memory of Professor Nicole Tomczak-Jaegermann.}

\begin{abstract}
In the last two decades a number of structures on the classical space of translation invariant valuations on convex bodies were discovered, e.g. product, convolution, a Fourier type transform.
In this paper a non-Archimedean analogue of the space of such (even) valuations with similar structures is constructed. It is shown that, like in the classical case, the new space equipped with either product
or convolution satisfies Poincar\'e duality and hard Lefschetz theorem.
\end{abstract}

\tableofcontents

\section{Introduction.}\label{S:introduction}
\begin{paragraphlist}
\item Theory of valuations on convex bodies (or convex valuations, for brevity) is a classical area of convex geometry, see e.g. \cite{alesker-kent} and \cite{schneider-book}, Ch. 6. By definition, valuations are finitely additive measures on the class of all convex compact subsets of a finite dimensional
{\itshape real} vector space.  In the last two decades a number of new structures on the space of smooth translation invariant convex valuations have been discovered.

The goal of this paper is to imitate this space of (even) valuations and the known structures on it in the new context of vector spaces over non-Archimedean local fields. Note that in this case we do not have
an interpretation of valuations as measures on something. We just construct a new space carrying structures with formally similar properties. Nevertheless most of the constructions in this paper are motivated by the usual valuations theory.

\item For these reasons let us say briefly a few words on the valuations theory on convex subsets of real vector spaces and structures on them. Here we emphasize mostly formal aspects of the theory relevant to this paper.
Let $V$ be a finite dimensional {\itshape real} vector space of dimension $n$. Let $Val^\infty(V)$ denote the space of smooth translation invariant valuations on convex compact subsets of $V$
(see e.g. \cite{alesker-04} or the lecture notes \cite{alesker-barcelona} for the definitions).  This space has a natural grading, called McMullen's decomposition
$$Val^\infty(V)=\oplus_{i=0}^n Val_i^\infty(V).$$
$Val_0^\infty(V)$ and $Val_n^\infty(V)$ are 1-dimensional, while all other summands are infinite dimensional.

The author  discovered a product on $Val^\infty(V)$ \cite{alesker-04} and a Fourier type transform (\cite{alesker-jdg-03}, \cite{alesker-ijm}).
In the recent preprint \cite{faifman-wannerer} Faifman and Wannerer  found a simpler approach to the Fourier type transform
on valuations.

Bernig and Fu \cite{bernig-fu-convol} discovered a convolution on convex valuations. The Fourier type transform intertwines product and convolution \cite{alesker-ijm}.

The algebra $Val^\infty(V)$ equipped either with product or convolution is a commutative associative graded algebra with a unit. The Fourier transform establishes an isomorphism between these two algebras.

These two algebras satisfy a hard Lefschetz type theorem. In some special cases it was proved by the author \cite{alesker-jdg-03}, \cite{alesker-gafa-HL}, \cite{alesker-ijm} and Bernig and Br\"ocker \cite{bernig-broecker}.
In the very recent preprint \cite{bernig-kotrbaty-wannerer} Bernig, Kotrbat\'y, and Wannerer proved a rather general hard Lefschetz type theorem in the language of convolution previously conjectured by Kotrbat\'y \cite{kotrbaty}.

Furthermore recently Kotrbat\'y \cite{kotrbaty} formulated a general conjecture on Hodge-Riemann bilinear relations for valuations in the language of convolution. He proved it in a special case.
In the above mentioned recent preprint \cite{bernig-kotrbaty-wannerer} by Bernig, Kotrbat\'y, and Wannerer this conjecture was fully proved; this preprint is based on another important special case of the conjecture
proved  by Kotrbat\'y and Wannerer \cite{kotrbaty-wannerer}.

\item Valuations and these structures on them found a number of non-trivial applications.
Theory of valuations has traditionally strong connections to integral geometry, see e.g. \cite{klain-rota}.
The above mentioned structures on convex valuations greatly enriched these connections and started to play a central role in various problems of integral geometry \cite{bernig-aig}, \cite{bernig-fu-annales}, \cite{bernig-fu-solanes-hermit}, \cite{fu-aig}.

The Hodge-Riemann bilinear relations in valuations formulated in the language of convolution imply the classical Alexandrov-Fenchel inequality for convex bodies  \cite{kotrbaty}. In some cases the Hodge-Riemann bilinear relations
can also be formulated in the language of product on valuations \cite{kotrbaty}, \cite{kotrbaty-wannerer-af}, and they imply new inequalities for mixed volumes \cite{alesker-ineq}, \cite{kotrbaty-wannerer-af}.

Some of the recent developments in valuations theory found deep connections to pseudo-Riemannian geometry  \cite{alesker-faifman}, \cite{bernig-faifman}, \cite{bernig-faifman-solanes2}, \cite{bernig-faifman-solanes}.

\item While valuations were originally introduced in convex geometry partly motivated by the needs of integral geometry, in the last two decades there were attempts to generalize the valuations theory beyond convexity.
Thus the space of valuations was introduced by the author on an arbitrary smooth manifold as the space of finitely additive measures of a special form on sufficiently 'nice' subsets of a manifold, see \cite{alesker-gafa-mfld} and references therein.
This space retains several properties of the classical space of valuations on convex sets, in particular  the product on valuations makes sense in this generality \cite{alesker-fu}.

In a different direction, the notion of a valuation on various classes of functions was introduced and investigated  \cite{colesanti-ludwig-mussnig-1}, \cite{colesanti-ludwig-mussnig-2}, \cite{knoerr}, \cite{ludwig-adv-geom-11}.

 This paper can be considered as another attempt to extend the valuations theory beyond convexity.  A description of its main results will be given in the next section.

\item {\bf Acknowledgements.} I am very grateful to J. Bernstein who suggested in 2001  the basic idea to study non-Archimedean analogues of the space of convex valuations,
 before all the structures on valuations were discovered.
The construction of the space of valuations
in this context is essentially contained in our joint paper \cite{alesker-bernstein}, Section 2, and it is used in the present paper. I thank also A. Aizenbud and D. Gourevitch for useful correspondence.

\item This version of the paper is different from the published one. The main results are the same. The main difference is the proofs of several properties of the  multiplicative structure on valuations.
The published version contained a mistake inherited from \cite{alesker-bernstein}, where it was claimed that the representation of the group $GL(V)$ in the space of sections of the line bundle $\cl_i^V$ over
the corresponding Grassmannian (see Section \ref{S:space of valuations}, paragraph \ref{item:m-l}) is not irreducible for $2\leq i\leq n-2$. However, the latter representation is irreducible.

\end{paragraphlist}


\section{Main results.}\label{S:main-results}
\begin{paragraphlist}
\item Let $\FF$ be a non-Archimedean local field. Let $V$ be an $n$-dimensional vector space over $\FF$ (see Section \ref{S:local-fields} for a reminder on definitions and basic properties). In Section \ref{S:space of valuations} we introduce the main objects of this paper: a complex vector space $Val^\infty(V)$
whose elements are called smooth valuations, and a complex vector space $Val(V)$ whose elements are called continuous valuations. The former space has no topology, while the latter is a Banach space.
Both spaces are infinite dimensional provided $n>1$, and $Val^\infty (V)$ is a dense subspace of $Val(V)$.
\begin{remark}
For the reader familiar with the classical theory of convex valuations, the notations $Val^\infty(V)$ and $Val(V)$ might be misleading. Both are analogous to even valuations rather than arbitrary ones.
Moreover the definition of $Val(V)$ is closer not to the definition of  continuous even convex valuations, but rather to Klain continuous even convex valuations in the sense of Bernig and Faifman \cite{bernig-faifman}.

Nevertheless we will keep this notation for the sake of simplicity.
\end{remark}

\item By definition, both above spaces are graded:
\begin{eqnarray*}
Val^\infty(V)=\oplus_{i=0}^n Val_i^\infty(V),
\end{eqnarray*}
and similarly
$$Val(V)=\oplus_{i=0}^n Val_i(V).$$
Here $Val_0^\infty(V)=Val_0(V)=\CC$, while $Val_n^\infty(V)=Val_n(V)$ is the 1-dimensional complex vector space of $\CC$-valued Lebesgue measures on $V$. All other spaces $Val_i^\infty(V), Val_i(V)$ are infinite dimensional (for $n>1$).

The space $Val_i^\infty(V)$ was first defined in \cite{alesker-bernstein}, Section 2, (see also Section \ref{S:space of valuations} below)
as the space of smooth sections
of certain $GL(V)$-equivariant complex line bundles over a Grassmannian over $\FF$.
It was shown in \cite{alesker-bernstein} that $Val_i^\infty(V)$ is an irreducible $GL(V)$-module. This definition is motivated by the analogy with the result from the same paper
\cite{alesker-bernstein}, Section 1, in the convex case.
Note that in the case of convex valuations the irreducibility of the space of valuations of given degree of homogeneity and parity (even or odd) was previously proved by the author \cite{alesker-irreduc}.

The space $Val_i(V)$ is defined in Section \ref{S:space of valuations} as the space of continuous sections of the same complex line bundle over a Grassmannian.

\item Let $X,Y$ be finite dimensional vector spaces over $\FF$. In Section \ref{S:Ext-product} we define the exterior product as a bilinear map
$$\boxtimes\colon Val(X)\times Val^\infty(Y)\to Val(X\times Y)$$
which is continuous with respect to the first argument. The exterior product on valuations on convex sets was defined in \cite{alesker-04} by the author.

Note that the exterior product of smooth valuations does not have to be smooth.
\begin{example}
1) Let $1\in Val_0(V)=\CC$. We will denote this element also by $\chi$ or $\chi_X$ to keep the analogy with the classical (convex) case where it corresponds to the Euler characteristic. Then
$$\chi_X\boxtimes \chi_Y=\chi_{X\times Y}.$$
\newline
2) Let $\mu_X,\mu_Y$ be Lebesgue measures on $X,Y$ respectively. Then $\mu_X\boxtimes\mu_Y$ is the product measure in the usual sense.
\end{example}

\item Let $F\colon X\to Y$ be a linear map of finite dimensional vector spaces over $\FF$. We construct the pull-back map
$$F^*\colon Val(Y)\to Val(X)$$
which is a continuous linear map of Banach spaces. The pull-back map on convex valuations was introduced by the author in \cite{alesker-ijm}.

\begin{example}
$F^*(\chi_Y)=\chi_X$.
\end{example}

If $F$ is injective then $F^*$ preserves the class of smooth valuations, see Theorem \ref{T:pull-back-thm}(4). In general it is not true.
The main properties of the pull-back are summarized in Theorem \ref{T:pull-back-thm}.

\item Let $D(V)=Val_n(V)$ denote the 1-dimensional complex vector space of $\CC$-valued Lebesgue measures on $V$, $n=\dim V$. In Section \ref{S:fourier-transform} we construct an analogue
of the Fourier transform
$$\FF\colon Val(V)\to Val(V^\vee)\otimes D(V),$$
where $V^\vee$ is the dual space of $V$. $\FF$ is an isomorphism of Banach spaces commuting with the action of $GL(V)$. It induces an isomorphism on the spaces of smooth valuations.
Our construction is a straightforward generalization of
the construction in \cite{alesker-jdg-03} for even convex valuations.

We show (Theorem \ref{T:fourier-ext-prod-commut}) that for smooth valuations $\phi,\psi$ one has
$$\FF\phi\boxtimes \FF\psi=\FF(\phi\boxtimes \psi).$$
For convex valuations this formula was recently proved in \cite{faifman-wannerer}.

\item For a linear map $F\colon X\to Y$ we define in Section \ref{S:push-forward} the push-forward map
$$F_*\colon Val(X)\otimes D(X)^*\to Val(Y)\otimes D(Y)^*$$
which is a continuous linear map of Banach spaces. By the definition
$$F_*=\FF\circ (F^\vee)^*\circ \FF^{-1},$$
where $F^\vee\colon Y^\vee\to X^\vee$ is the dual map. For convex valuations the push-forward was defined in \cite{alesker-ijm}.

\item In Section \ref{S:product} we define product on the space of smooth valuations $Val^\infty(V)$ as follows
$$\phi\cdot \psi:=\Delta^*(\phi\boxtimes\psi),$$
where $\Delta\colon V \inj V\times V$ is the diagonal imbedding, i.e. $\Delta(v)=(v,v)$. It is shown that equipped with this product, $Val^\infty(V)$ is a commutative associative
algebra with a unit (equal to $\chi_V$).  It is graded:
$$Val_i^\infty(V)\cdot Val_j^\infty(V)\subset Val_{i+j}^\infty(V).$$
Denote $n:=\dim V$. $Val^\infty(V)$ satisfies Poincar\'e duality: the bilinear map given by the product
$$Val_i^\infty(V)\times Val_{n-i}^\infty(V)\to Val^\infty_n(V)=D(V)$$
is a perfect pairing, i.e. for any $0\ne \phi \in Val^\infty_i(V)$ there exists $\psi\in Val_{n-i}^\infty(V)$ such that $\phi\cdot \psi\ne 0$.

The product on convex valuations was introduced by the author \cite{alesker-04}.

\item In Section \ref{S:hard-Lefschetz} it is shown that the algebra of smooth valuations $Val^\infty(V)$ satisfies a version of hard Lefschetz theorem.
To state it, let us denote by $V_1\in Val_1^\infty(V)$ the only (up to a proportionality) element invariant under a maximal compact subgroup $GL_n(\co)$ of $GL_n(\FF)\simeq GL(V)$.
Let $0\leq i<n/2$. Then the map $Val_i^\infty(V)\to Val_{n-i}^\infty(V)$ given by
$$\phi\mapsto \phi\cdot (V_1)^{n-2i}$$
is an isomorphism.

The proof of this theorem uses properties of the Radon transform on Grassmannians over $\FF$ due to Petrov and Chernov \cite{petrov-chernov}.

\item In Section \ref{S:convolution} we introduce a convolution
$$\ast\colon (Val^\infty(V)\otimes D(V)^*)\times (Val^\infty(V)\otimes D(V)^*)\to Val^\infty(V)\otimes D(V)^*$$
by $\phi\ast\psi=a_*(\phi\boxtimes \psi)$, where $a\colon V\times V\to V$ is the addition map, i.e. $a(x,y)=x+y$.

By Proposition \ref{P:convol-vs-product} the convolution is related to the product and the Fourier transform by the formula
$$\FF\phi\ast\FF\psi=\FF(\phi\cdot\psi).$$

Convolution also satisfies Poincar\'e duality and hard Lefschetz type theorem  (Theorem \ref{T:convolution-properties}).

On convex valuations the convolution was introduced by Bernig and Fu \cite{bernig-fu-convol}.

\item An interesting open question is to establish for $Val^\infty(V)$ the Hodge-Riemann bilinear
relations similar to \cite{kotrbaty} (see also \cite{kotrbaty-wannerer}, \cite{bernig-kotrbaty-wannerer}).

\end{paragraphlist}

\section{Reminder on local fields.}\label{S:local-fields}
\begin{paragraphlist}
\item In this section we collect a few basic well known facts on local fields sufficient for this paper. We refer to \cite{weil}, Ch. 1, for details.

By definition, a local field is a topological locally compact non-discrete field. There is a classification of such fields: they are precisely $\RR,\CC, \FF_q((t))$, and finite extensions of the
fields of $p$-adic numbers $\QQ_p$. Here $\FF_q$ denotes the finite field with $q$ elements, and $\FF_q((t))$ denotes the field of formal Laurent power series.
The first two examples, namely $\RR,\CC$, are called Archimedean, while all others are called non-Archimedean local fields.

\item Let $\FF$ be a non-Archimedean local field. It has a unique maximal compact subring $\co\subset \FF$. For example if $\FF=\FF_q((t))$ then $\co=\FF_q[[t]]$ is the ring of all Taylor power series.
If $\FF=\QQ_p$ then $\co=\ZZ_p$ is the ring of $p$-adic integers.

The field of fractions of $\co$ equals $\FF$.

\item $\co$ has a unique maximal ideal $\frak{m}\subset \co$. For example for $\FF=\FF_p((t))$ the ideal $\frak{m}$ is generated by $t$,
while for $\FF=\QQ_p$ the ideal $\frak{m}$ is generated by $p$.

The quotient $k:=\co/\frak{m}$ is necessarily a finite field; it is called the residue field of $\FF$.

\item There exists a unique multiplicative norm
$$ |\cdot|\colon \FF\to \RR_{\geq 0}$$
such that
\begin{eqnarray*}
|x|=1 \,\,\, \forall x\in\co\backslash \frak{m},\\
|x|=\frac{1}{|k|^{i}} \,\,\, \forall x\in \frak{m}^i\backslash \frak{m}^{i+1}, \mbox{where } i\geq 1.
\end{eqnarray*}
where $|k|$ denotes the cardinality of the residue field.
Multiplicativity means that  $|x\cdot y| =|x|\cdot |y|$ for any $x,y\in \FF$.

This norm satisfies the strengthened triangle inequality
$$|x+y|\leq \max\{|x|,|y|\}.$$

\item The norm $|\cdot|$ has the following property. Let $\mu$ be a Lebesgue measure on $\FF$ ($\mu$ exists and is unique up to a proportionality). Let $x\in \FF$.
Then $$\mu(x\cdot A)=|x|\mu(A)$$
for any compact subset $A\subset \FF$.

\end{paragraphlist}

\section{Lattices over non-Archimedean local fields.}\label{S:lattice}
\begin{paragraphlist}
\item In this section $\FF$ denotes a non-Archimedean local field, and $\co\subset \FF$ its ring of integers. In this section we review, mostly following \cite{weil}, a few well known facts on finite dimensional $\FF$-vector spaces
and lattices in them.

A proof of the following result can be found in \cite{schaeffer}, Thm. 3.2, Ch. 1.
\begin{theorem}\label{T:isomorphism-finite-dim}
Let $V$ be an $n$-dimensional Hausdorff topological vector space over the local field $\FF$. Let $v_1,\dots,v_n$ be its basis. Then the map $\FF^n\to V$ given by
$$(x_1,\dots,x_n)\mapsto x_1v_1+\dots+x_nv_n$$
is an isomorphism of topological vector spaces when the source space is equipped with the product topology.
\end{theorem}

\item Let $V$ be an $n$-dimensional Hausdorff topological vector space over the local field $\FF$.
\begin{definition}\label{D:lattice}
A lattice $L$ in $V$ is a compact open $\co$-submodule of $V$.
\end{definition}

 \begin{lemma}\label{L:sublattice-quotient}
Let $L\subset V$ be a lattice. Let $E\subset V$ be a vector subspace. Then $E\cap L$ is a lattice in $E$, and $L/E\cap L$ is a lattice in $V/E$.
\end{lemma}
{\bf Proof.} This immediately follows from Definition \ref{D:lattice}. \qed

\begin{theorem}[\cite{weil}, Ch. II, \S 2, Thm. 1]\label{T:free-lattice}
Let $L\subset V$ be a lattice.
\newline
(i) Then $V$ has a basis $v_1,\dots, v_n$ such that $L=\co v_1\oplus\dots \oplus \co v_n$. In particular $L$ is a free module of rank $n$.
\newline
(ii) Moreover if $\{0\}=V_0\subset V_1\subset \dots\subset V_{n-1}\subset V_n=V$ is a sequence of linear subspaces such that $\dim V_i=i$. Then the above vectors $ v_1,\dots, v_n$
can be chosen so that $v_{1},\dots,v_i$ is a basis of $V_i$ for any $i$.
\end{theorem}

\begin{remark}\label{Rem:basis-of-lattice}
In the assumptions of part (ii) of the last theorem, one clearly has for each $i$
$$L\cap V_i=\co v_1\oplus\dots\oplus \co v_i.$$
\end{remark}

\item Given a lattice $L\subset V$. Denote by $GL(L)$ the subgroup
$$GL(L):=\{T\in GL(V)|\, T(L)=L\}.$$

\begin{proposition}\label{P:transitivity}
Let $L\subset V$ be a lattice. The natural action of the group $GL(L)\simeq GL_n(\co)$ on the Grassmannian $Gr_i^V$ is transitive.
\end{proposition}
{\bf Proof.} We may and will assume that $V=\FF^n$, $L=\co^n$. Let $e_1,\dots, e_n\in \FF^n$ be the standard basis. Let $E_0:=span\{e_1,\dots,e_i\}\in Gr_i^V$.
Let $E\in Gr_i^V$. We have to show that there exists $T\in GL_n(\co)$ such that $T(E_0)=E$.

By Lemma \ref{L:sublattice-quotient} and Theorem \ref{T:free-lattice} there exists a basis $v_1,\dots,v_i$ of $E$ such that $$E\cap \co^n=\co v_1\oplus\dots\oplus \co v_i.$$
Similarly there exists a basis $\bar v_{i+1},\dots, \bar v_{n}$ of $V/E$ such that
$$L/L\cap E=\co \bar v_{i+1}\oplus \dots\oplus \co \bar v_{n}.$$
Let us choose $v_j\in L$, $j=i+1,\dots,n$, such that $v_j\equiv \bar v_j mod (E)$.
It is easy to see that $\co^n=\co v_1\oplus \dots\oplus\co v_n$.

Then define $T\colon \FF^n\to \FF^n$ by $T(e_j)=v_j$, $j=1,\dots,n$. Clearly $T(E_0)=E$ and $T(\co^n)=\co^n$, i.e. $T\in GL_n(\co)$. \qed

\item \begin{proposition}\label{P:close-norms}
Let $L\subset V$ be a lattice. Let $F_0\colon W\inj V$ be an injective linear map of vector spaces. There exists a neighborhood $U$ of $F_0$ in $Hom(W,V)$ such that any $F\in U$ is also injective
and $$F^{-1}(Im(F)\cap L)=F_0^{-1}(Im(F_0)\cap L).$$
\end{proposition}

{\bf Proof.} We may and will assume that $V=\FF^n$ and $L=\co^n$. Let us define a norm on $\FF^n$ by
$$||(x_1,\dots,x_n)||=\max_{i}|x_i|.$$
Clearly
\begin{eqnarray*}
||x+y||\leq \max\{||x||,||y||\},\\
||\lam\cdot x||=|\lam|\cdot ||x||.
\end{eqnarray*}
Then $$\co^n=\{x\in \FF^n|\,\,\, ||x||\leq 1\}$$
is the unit ball of this norm.

By Theorem \ref{T:free-lattice}(ii) we may also assume that
$$Im(F_0)=span\{e_1,\dots,e_i\},$$
where $e_1,\dots,e_n\in \FF^n$ is the standard basis.
Let $w_1,\dots,w_i\in W$ be the basis such that $F_0(w_j)=e_j$ for $1\leq j\leq i$.

Then clearly
$$F_0^{-1}(Im(F_0)\cap L) =\{\sum_{j=1}^i x_jw_j|\,\, \max_{1\leq j\leq i}|x_j|\leq 1\}.$$

Define $U:=\{F\in Hom(W,\FF^n)|\,\,\, ||F(w_j)-e_j||<1\,\, \forall 1\leq j\leq i\}$.
Let $F\in U$. The vectors $F(w_1),\dots,F(w_i)$ are linearly independent since these are vectors with coordinates from $\co$, and their reduction modulo the maximal ideal $\mathfrak{m}$ of $\co$
are the first $i$ vectors of the standard basis of $(\co/\mathfrak{m})^n$. It follows that $F$ is injective.

Now it remains to show that
$$F^{-1}(Im(F)\cap L)=\{\sum_{j=1}^ix_jw_j|\,\max_{1\leq j\leq i}|x_j|\leq 1\}.$$
Equivalently, one has to show that for $x_1,\dots,x_i\in \FF$ the inequality
\begin{eqnarray}\label{E:inequal-ball}
||F(x_1w_1+\dots+x_iw_i)||\leq 1
\end{eqnarray}
holds if and only if $\max_{1\leq j\leq i}\{|x_j|\}\leq 1$.

The 'if' part follows since $||Fw_j||\leq \max\{||Fw_j-e_j||,||e_j||\}=1$ for any $1\leq j\leq i$.

Conversely, let us assume that (\ref{E:inequal-ball}) holds. Without loss of generality we may assume that $|x_1|=\max_{1\leq j\leq i}\{|x_j|\}$.
Let us denote $\theta_j=F(w_j)-e_j$, $1\leq j\leq i$. Then $||\theta_j||<1$. We have
\begin{eqnarray*}
1\geq |x_1|\cdot ||(e_1+\theta_1)+\sum_{j=2}^i\frac{x_j}{x_1}(e_j+\theta_j)||=\\
|x_1|\cdot ||(1,\frac{x_2}{x_1},\dots,\frac{x_i}{x_1},0,\dots,0)+(\theta_1+\sum_{j=2}^i\frac{x_j}{x_1}\theta_j)||.
\end{eqnarray*}
Since $|\frac{x_j}{x_1}|\leq 1$ the norm of the first summand equals 1, while the norm of the second summand (in the parenthesis) is strictly less than 1.
Hence the norm of their sum is equal to 1. Thus we get $1\geq |x_1|\cdot 1=|x_1|$. \qed

\end{paragraphlist}

\section{Lebesgue measures on vector spaces.}\label{S:measures}
\begin{paragraphlist}
\item\label{item-Leb-measures-1} In this section we assume that all vector spaces are finite dimensional over a non-Archimedean local field $\FF$.
For such a vector space $V$ we denote by $D(V)$ throughout the article the (one-dimensional) complex vector space of $\CC$-valued Lebesgue measures on $V$.

Let $\sigma\colon V\to W$ be a surjective linear map  between such vector spaces. Let $K=Ker(\sigma)$. We are going to construct a linear isomorphism
\begin{eqnarray}\label{E:iso-measures}
\tilde\sigma\colon D(K)\otimes D(W)\tilde\to D(V).
\end{eqnarray}
Let $\mu_K\in D(K), \mu_W\in D(W)$. Since measures are linear functionals on compactly supported continuous functions, define for any $\phi\in C_c(V)$
\begin{eqnarray}\label{E:iso-measures-2}
\int \phi(v)  d\tilde\sigma(\mu_K\otimes \mu_W):=\int_{w\in W} d\mu_W(w)\int_{k\in \sigma^{-1}(w)}\phi(k) d\mu_K(k),
\end{eqnarray}
where we identify the measure $\mu_K$ on $K$ with its (arbitrary) translate to the parallel affine subspace $\sigma^{-1}(w)$.

An equivalent description is as follows. Let us choose a splitting $V=K\oplus L$. Then
$$\sigma|_L\colon L\to W$$
is an isomorphism. Set $\mu_L:=(\sigma|_L^{-1})_*\mu_W$. Then it is easy to see that
$$\tilde\sigma(\mu_K\otimes \mu_W)=\mu_K\boxtimes \mu_L,$$
where $\boxtimes$ denotes the usual product measure.

\item Let $X^\vee$ denote the dual space of a vector space $X$. The goal of this paragraph is to construct a canonical isomorphism
\begin{eqnarray}\label{E:isom-canon-dens}
D(X)^*\tilde\to D(X^\vee).
\end{eqnarray}
Since $$Hom_\CC(D(X)^*,D(X^\vee))\simeq D(X)\otimes D(X^\vee),$$ it suffices to construct a canonical non-zero element in $D(X)\otimes D(X^\vee)$.

For a lattice $\Lam\subset X$ let us define the dual lattice in $X^\vee$
$$ \Lam^\vee:=\{f\in X^\vee|\, f(\Lam)\subset \co\}.$$
It is easy to see that $\Lam^\vee$ is a lattice in $X^\vee$.
\begin{lemma}\label{L:dual-measure-non-arch}
Let $\mu\in D(X)$ be a non-vanishing ($\CC$-valued) Lebesgue measure on $X$. Then there is a unique Lebesgue measure on $X^\vee$ denoted by $\mu^{-1}$
such that for any lattice $\Lam\subset X$ one has
\begin{eqnarray}\label{E:measures-equal}
\mu^{-1}(\Lam^\vee)=\frac{1}{\mu(\Lam)}.
\end{eqnarray}
\end{lemma}
{\bf Proof.} Let us fix a lattice $\Lam\subset X$. We can obviously construct a unique $\mu^{-1}$ such that  (\ref{E:measures-equal}) is satisfied for this $\Lam$.
It remains to show that  (\ref{E:measures-equal}) is satisfied for any other lattice $\tilde\Lam$.
For there exists $g\in GL(X)$ such that $\tilde\Lam=g(\Lam)$. Then $\tilde\Lam^\vee=(g^\vee)^{-1}\Lam^\vee$, where $g^\vee$ is the dual map of $g$.
Then we have
\begin{eqnarray*}
\mu^{-1}(\tilde\Lam^\vee)=\mu^{-1}((g^\vee)^{-1}\Lam^\vee)=\\
|\det g^\vee|^{-1}\mu^{-1}(\Lam^\vee)=\frac{1}{|\det g|\mu(\Lam)}=
\frac{1}{\mu(\tilde\Lam)}
\end{eqnarray*}
\qed

\hfill

Let us now construct the promised non-zero element of $D(X)\otimes D(X^\vee)$. Fix an arbitrary $\CC$-valued non-vanishing Lebesgue measure $\mu\in D(X)$.
Define the element to be
$$\mu\otimes \mu^{-1}\in D(X)\otimes D(X^\vee) .$$
The following claim is now obvious.
\begin{claim}\label{Cl:canon-element}
 The element $\mu\otimes \mu^{-1}$ is independent of $\mu$.
 It is $GL(X)$-invariant under the natural action of this group on $D(X)\otimes D(X^\vee)$.
\end{claim}
\end{paragraphlist}

Now we can explicitly describe the isomorphism (\ref{E:isom-canon-dens}). Fix a non-vanishing $\mu\in D(X)$. There is a unique element $\mu^\vee\in D(X)^*$ such that $\mu^\vee(\mu)=1$.
The isomorphism (\ref{E:isom-canon-dens}) maps $\mu^\vee\mapsto \mu^{-1}$. It is easy to see that this map is independent of $\mu$.

Very often, by the abuse of notation, we will write in this paper $$D(X)^*=D(X^\vee)$$ meaning the isomorphism (\ref{E:isom-canon-dens}).


\section{Reminder on analytic manifolds over local fields.}\label{S:analytic}
\begin{paragraphlist}
\item The goal of this section is to review very briefly the notion of analytic manifold over a non-Archimedean local field $\FF$. For more details we refer to \cite{serre}, part II.

The theory of analytic manifolds over such a field is similar to the theory of real analytic manifolds, at least at the basic level needed for this paper.
The main examples of analytic manifolds to keep in mind for the purposes of this paper are the Grassmannians of linear $i$-dimensional subspaces in $\FF^n$ and, more generally, the manifolds of partial flags in $\FF^n$.
All the material of this section is well known.

\item Let $U\subset \FF^n$ be an open subset. A function $f\colon U\to \FF$ is called analytic if any point $a\in U$ has a ball centered at $a$ in which $f$ can be represented by a power
 series which absolutely converges in this ball.

Let $F=(F_1,\dots, F_m)$ be a map $F\colon U\subset \FF^n\to \FF^m$.  $F$ is called analytic if every $F_i$ is analytic.

Composition of analytic maps is analytic (\cite{serre}, part II, Ch. II, Theorem 2).

For an analytic map $F$ its Jacobian is defined as usual
$$J(F)=\left(\frac{\pt F_i}{\pt x_j}\right).$$
Thus $J(F)$ is an $m\times n$ matrix whose entries are analytic functions.

\item\label{item:implicit-function} For an analytic map $$F\colon U\subset \FF^n\to \FF^n$$
a version of the inverse function theorem holds (see \cite{serre}, part II, Ch. III, \S 9). Namely assume that at a point $a\in U$
$$\det(J(F)_a)\ne 0.$$
Then there exists an open neighborhood $V$ of $a$ such that $F(V)\subset \FF^n$ is open, $F|_V\colon V\to F(V)$ is a homeomorphism, and $F^{-1}\colon F(V)\to V$ is an analytic map.

It is easy to see that the inverse function theorem implies in the usual way the implicit function theorem for analytic functions. It is formulated in the usual way and we leave it
to the reader. The implicit function theorem will be used in this paper in the proof of Lemma \ref{L:product-by-intrin-volume}.

\item A Hausdorff topological space $X$ is an analytic manifold if it admits an open covering $\{U_\alp\}$, homeomorphisms
$$\phi_\alp\colon U_\alp\to V_\alp,$$
where $V_\alp\subset \FF^n$ are open subsets such that the transition maps $$\phi_\alp|_{U\alp\cap U_\beta}\circ \phi_\beta^{-1}|_{\phi_\beta(U\alp\cap U_\beta)}$$ are analytic for any $\alp,\beta$.

One can define analytic maps between analytic manifolds in the obvious way.

\item\label{item:Lie-grp} A topological group $G$ which is also an analytic manifold is called a Lie group if
\newline
(a) the product map $G\times G\to G$ given by $(x,y)\mapsto x\cdot y$ is analytic;
\newline
(b) the inverse map $G\to G$ given by $x\mapsto x^{-1}$ is analytic.

Basic examples of Lie groups are $GL_n(\FF),GL_n(\co)$. Another example is the subgroup of $GL_n(\FF)$ stabilizing the given partial flag:
$$G=\left\{\left[\begin{array}{c|c|c|c}
             A_1&\ast&\dots&\ast\\
             \hline
             0&A_2&\dots&\ast\\
             \hline
             0&0&\ddots&\ast\\
             \hline
             0&0&\dots&A_s
             \end{array}\right]\right\},$$
             where $A_1,A_2,\dots, A_s$ are invertible square matrices. Thus $G$ is the subgroup of invertible block upper triangular matrices.

\item Let $G$ be a Lie group. Let $H\subset G$ be a subgroup which is an analytic submanifold (this notion is naturally defined). Then $H$ is also a Lie group which is called a Lie subgroup of $G$
(see \cite{serre}, part II, Ch. IV, \S 2.3).

\begin{theorem}[\cite{serre}, part II, Ch. IV, \S 5]\label{T:serre}
Let $G$ be a Lie group. Let $H\subset G$ be its Lie subgroup. Then $G/H$ has a unique structure of analytic manifold such that the natural map $G\to G/H$ is analytic and has surjective differential at every point.
\end{theorem}

This theorem implies immediately that the Grassmannians and, more generally, partial flag spaces are analytic manifolds.

\item Let us discuss now integration over analytic manifolds over a non-Archimedean local field $\FF$. Let us start with integration in $\FF^n$.

Let $dx$ denote the Lebesgue measure on $\FF^n$ normalized so that
its value on $\co^n$ equals 1.  Let $U,V\subset \FF^n$ be compact open subsets. Let
$$f\colon V\to\CC$$
be a continuous function. Let $F\colon U\tilde\to V$ be an analytic homeomorphism such that $F^{-1}$ is also analytic. Then there is the following change of variables formula (see \cite{igusa}, \S 7.4)
\begin{eqnarray}\label{E:change-formula}
\int_V f(y)dy=\int_U f(F(x))|\det J(F)_x|dx.
\end{eqnarray}

\item Let us define the complex line bundle $|\ome_X|$ over $X$ called the line bundle of densities. We will see that  its continuous sections can be integrated over $X$.
Fix an atlas of charts $\{(U_\alp,\phi_\alp)\}$ on $X$. Consider the transition functions
$$F_{\alp\beta} := \phi_\alp|_{U_\alp\cap U_\beta}\circ \phi^{-1}_\beta|_{\phi_\beta(U_\alp\cap U_\beta)}\colon \phi_\beta(U_\alp\cap U_\beta)\tilde\to  \phi_\alp(U_\alp\cap U_\beta).$$
Choose the trivial line bundle $U_\alp\times \CC$ over each $U_\alp$ and identify them over pairwise intersections $U_\alp\cap U_\beta$ as follows
$$(x,z)\sim (x,|\det J(F_{\alp\beta})_{\phi_\beta(x)}|\cdot z).$$

\item\label{item:integration} Let us assume in addition that $X$ is compact.
Then continuous sections of $|\ome_X|$ can be integrated over $X$ by patching together local integrations over subsets of $\FF^n$. More precisely let us fix a finite atlas of charts $\{(U_\alp,\phi_\alp)\}$ on $X$.
There exists a subordinate partition of unity $\{\psi_\alp\}$ (see e.g. \cite{rudin-real-complex}, Theorem  2.13), namely there exist continuous functions  $\psi_\alp\colon X\to \RR_{\geq 0}$  with $supp(\psi_\alp)\subset U_\alp$ and such that
$$\sum_\alp\psi_\alp=1.$$
Let $\ome$ be a continuous section over $X$ of the line bundle $|\ome_X|$. Then
$$\ome=\sum_\alp \psi_\alp\cdot \ome.$$
Clearly $supp( \psi_\alp\cdot \ome)\subset U_\alp$. Then $\int_{U_\alp}\psi_\alp\cdot \ome$ is well defined since $(U_\alp,\phi_\alp)$ is a chart. Then one defines
$$\int_X\ome:=\sum_\alp \int_{U_\alp}\psi_\alp \cdot \ome.$$
This number is independent of the atlas of charts and the subordinate partition of unity.

\end{paragraphlist}

\section{Representations of the group $GL_n(\co)$.}\label{S:representations}
\begin{paragraphlist}

\item In this section we summarize a few known results on representations of the group $GL_n(\co)$
in the space of complex valued functions on Grassmannians $Gr_k^{\FF^n}$. No result of this section is novel.

Let $\FF$ be a non-Archimedean local field. Let $\co\subset \FF$ be its ring of integers. The group $GL_n(\co)$ is compact and acts transitively on the Grassmannian $Gr_k^{\FF^n}$. Hence there is a unique probability (Haar) measure $\mu_{Haar}$
on $Gr_k^{\FF^n}$ invariant under this group. By the general representation theory of compact groups the representation of $GL_n(\co)$ in $L^2(Gr_k^{\FF^n})$ is unitary and has a dense subspace which is an orthogonal countable  direct sum of irreducible representations.
The irreducible representations are necessarily finite dimensional.

\item\label{item:mult-one} Each irreducible representation of $GL_n(\co)$ enters $L^2(Gr_k^{\FF^n})$ with multiplicity at most 1. For $char(\FF)=0$ this was proven first in \cite{hill-94}, Corollary 3.2,  in general in \cite{bader-onn-2012}.

\item The linear subspace $C^\infty(Gr_k^{\FF^n})$ of locally constant (called smooth) functions is $GL_n(\co)$-invariant and dense in $L^2(Gr_k^{\FF^n})$.

Every finite dimensional $GL_n(\co)$-invariant subspace of $L^2(Gr_k^{\FF^n})$
is contained in $C^\infty(Gr_k^{\FF^n})$; moreover the representation of $GL_n(\co)$ in this subspace factors through the quotient of $GL_n(\co)$ by an open normal subgroup of finite index.
 The last two statements follow from the fact that, as a topological group, $GL_n(\co)$ is a pro-finite
group (i.e. inverse limit of finite groups).

\end{paragraphlist}

\section{Space of valuations.}\label{S:space of valuations}

\begin{paragraphlist}
\item Let $V$ be an $n$-dimensional vector space over a non-Archimedean local field $\FF$. We introduce the main object of study $Val(V)=\oplus_{k=0}^n Val_k(V)$ which is a vector space over $\CC$ and is an analogue of
{\itshape even and Klain continuous} valuations on convex sets in the terminology of \cite{bernig-faifman}, Section 3.

We define $\Vai_0(V):=\CC$ and $\Vai_n(V):=D(V)$, the latter is the (1-dimensional) space of complex valued Lebesgue measures on $V$. Let now $1\leq k\leq n-1$.
We are going to define in this section the space $\Vai_k(V)$.

\item\label{it:CT} Let us start with an elementary construction. Given two vector spaces $X$ and $Y$ of equal (finite)  dimension over the local field $\FF$. Define a map linear with respect to the second variable
\begin{eqnarray}\label{E:pull-back}
\ct\colon Hom_\FF(X,Y)\times D(Y)\to D(X).
\end{eqnarray}
Let $F\in Hom_\FF(X,Y), \mu\in D(Y)$.  There exists a unique Lebesgue measure $\nu$ on $X$ such that for some (equivalently, any) compact subset $A\subset X$ with non-empty interior one has
$\nu(A)=\mu(F(A))$. Define $\ct(F,\mu):=\nu$.

\begin{claim}\label{E:properties-T}
(1) The map $\ct$ is linear with respect to the second argument.
\newline
(2) If $F$ is invertible then $\ct(F,\mu)=(F^{-1})_*(\mu)$, where $G_*$ denotes the push-forward of measures under a map $G$. Otherwise $\ct(F,\mu)=0$.
\newline
(3) The map $\ct$ is jointly continuous.
\end{claim}
A proof is left to the reader.

\item\label{item:m-l} Let $\cl_k^V\to Gr_k^V$ be the complex line bundle whose fiber over $E\in Gr_k^V$ is equal to the 1-dimensional space of complex valued Lebesgue measures on $E$.
\begin{definition}
Define the space of smooth $k$-homogeneous valuations $$Val^\infty_k(V)=C^\infty(Gr_k^V,\cl_k^V),$$ and the space of continuous $k$-homogeneous valuations
$$Val_k(V)=C(Gr_k^V,\cl_k^V).$$
\end{definition}

$Val^\infty_k(V)$ is (algebraically) irreducible $GL(V)$-module, i.e. has no proper $GL(V)$-invariant subspaces. $Val_k(V)$ is topologically irreducible $GL(V)$-module,
i.e. has no proper closed $GL(V)$-invariant subspaces.

Let $|\ome_{n-k}^V|\to Gr_{n-k}^V$ denote the line bundle of densities on $Gr_{n-k}^V$ as defined in Section \ref{S:analytic}.
Its continuous global sections can be integrated over the Grassmannian $Gr_{n-k}^V$.

Let us denote by $\cm'_{n-k}\to Gr_{n-k}^V$ the line bundle whose fiber over $F\in Gr_{n-k}^V$ is equal to the space of complex valued Lebesgue measures on $V/F$.
Set finally $\cm_{n-k}^V:=\cm'_{n-k}\otimes |\ome_{n-k}^V|$.
All the line bundles $\cl_{n-k}^V,|\ome_{n-k}^V|,\cm'_{n-k},\cm_{n-k}^V$ are $GL(V)$-equivariant in a natural way.

A section of any of the above vector bundles is called smooth if its stabilizer in $GL(V)$ is an open subgroup. More generally a vector in a continuous representation of
$GL(V)$ is called smooth if its stabilizer is open. It is easy to see that any function on $Gr_i^V$ is smooth in this sense if and only if it is locally constant.

\item\label{item-cD} Let us define the $GL(V)$-equivariant operator between spaces of smooth sections
\begin{eqnarray}\label{E:operator-D}
\cd\colon C^\infty(Gr_{n-k}^V,\cm_{n-k}^V)\to C^\infty(Gr_k^V,\cl_k^V)
\end{eqnarray}
as follows. For $F\in Gr_{n-k}^V$ and $E\in Gr_k^V$ define $p_{E,F}\colon E\to V/F$ the natural map.
 Given a section $\xi\in C^\infty(Gr_{n-k}^V,\cm_{n-k}^V)$, define for any $E\in Gr_k^V$
\begin{eqnarray}\label{E:def-cD}
\cd(\xi)(E)=\int_{F\in Gr_{n-k}^V}\ct(p_{E,F},\xi(F)),
 \end{eqnarray}
 where $\ct$ is the map from paragraph \ref{it:CT}. Note that $$\ct(p_{E,F},\xi(F))\in D(E)\otimes |\ome_{n-k}^V|\big|_F.$$ Hence the integral belongs to $D(E)$.

Let us reformulate the definition of $\cd$. Given a $k$-dimensional subspace $E$.
The set of $F\in Gr_{n-k}^V$ intersecting $E$ non-transversally, or equivalently $p_{E,F}$ is non-invertible, has zero measure in the Grassmannian. Ignoring this subset we have
\begin{eqnarray}\label{E:iii1}
\cd(\xi)(E)=\int_{F\in Gr_{n-k}^V}(p_{E,F})^{-1}_*\xi(F).
\end{eqnarray}

\begin{proposition}\label{P:irred}
(1) The operator $\cd$ is $GL(V)$-equivariant. It extends uniquely by continuity to the space of continuous sections. Then it maps continuous sections to continuous, smooth to smooth.
\newline
(2) The source and the target of $\cd$ are irreducible $GL(V)$-representations. Consequently, $\cd$ is an isomorphism between them.
\end{proposition}

Part (1) follows from Claim \ref{E:properties-T}(3). Part (2) was proved in \cite{alesker-bernstein}, Theorem 2.1; see also \cite{gourevitch-compos}, Corollary 1.3, for a more general statement.

\item Let us fix a lattice $\Lam\subset V$. Since the action of $GL(\Lam)$ on $Gr_k^{V}$ is transitive, there exist  unique (up to a proportionality) non-zero $GL(\Lam)$-invariant continuous sections of $\cm_{n-k}^V$ and of $\cl_k^V$ (which are obviously smooth).
It is easy to see that $\cd$ applied to the former is a non-zero multiple of the latter. We call such a $GL(\Lam)$-invariant section a spherical vector.

\end{paragraphlist}

\section{Pull-back on valuations.}\label{S:pull-back}
\begin{paragraphlist}
\item Let $\FF$ be a non-Archimedean local field. The goal of this section is to construct an operation of pull-back on valuations. More precisely we will prove
\begin{theorem}\label{T:pull-back-thm}
Let $X$ and $Y$ be finite dimensional vector spaces over $\FF$. For any linear map $F\colon X\to Y$ there exists a canonical continuous linear map
$$F^*\colon Val(Y)\to Val(X)$$
which is called the pull-back map and satisfies the following properties:
\newline
(1) $F^*$ preserves degree of homogeneity, i.e. $F^*(Val_k(Y))\subset Val_k(X)$;
\newline
(2) $(F\circ G)^*=G^*\circ F^*$;
\newline
(3) $Id^*=Id$;
\newline
(4) If $F$ is injective then $F^*$ preserves the subspace of smooth valuations, i.e. $F^*(Val^\infty(Y))\subset Val^\infty(X)$.
\end{theorem}

Proof of this theorem occupies the rest of this section.

\item The map $F^*$ on valuations of degree 0 is just the identity map of $\CC$.

\item Let $F\colon X\to Y$ be a linear map and $\dim X=\dim Y=n$. The pullback $F^*\colon Val_n(Y)\to Val_n(X)$ is the map $F^*\colon D(Y)\to D(X)$ given, by definition, by
$$F^*(\mu)=\ct(F,\mu),$$
where $\ct$ is the map from Claim \ref{E:properties-T}.

\item Let $F\colon X^n\to Y^m$ be a linear map. Let us define first the pull-back map
$$F^*\colon D(Y)\to Val_m(X).$$
If $m>n$ then $F^*=0$ by the definition. Assume that $m\leq n$.
For $\mu\in D(Y)$ and any $E\in Gr_m^X$ set
$$(F^*\mu)(E):=\ct(F|_E,\mu).$$
The continuity of $\ct$ with respect to the first variable implies that $F^*\mu$ is continuous.
\item \begin{lemma}\label{L:pull-back-density}
The linear span of valuations of the  form $F^*\mu$ over all possible linear maps $F\colon X^n\to Y^m$ and $\mu\in D(Y)$ is dense in $Val_{m}(X)$.
\end{lemma}
{\bf Proof.}
It is clear that the linear span of valuations of the form $F^*\mu$ is a $GL(X)$-invariant subspace of $Val_{m}(X)$.
Since the latter space is $GL(X)$-irreducible by Proposition \ref{P:irred}, the result follows.
\qed

\hfill

The constructed map $F^*\colon D(Y)\to Val_{\dim Y}(X)$ is called pull back on densities.

\item The following corollary will be used later.
\begin{corollary}\label{COR:constr-valuat}
Let $\cx$ be a compact metrizable space. Let $m$ be a complex valued Borel measure on $\cx$. Let
$$T\colon \cx\to Hom(X,Y)$$
be a continuous map when $X,Y$ are finite dimensional vector spaces over $\FF$. Let $\mu \in D(Y)$. Then
$$\int_\cx [T(x)]^*(\mu) dm(x)$$
belongs to $Val_{\dim Y}(X)$.
\end{corollary}
{\bf Proof.} By Claim \ref{E:properties-T} the map $\cx\to C(Gr_{\dim Y}^X,\cl_{\dim Y}^X)$ given by the expression under the last integral $x\mapsto [T(x)]^*(\mu)$ is continuous.
By Lemma \ref{L:pull-back-density} the expression under the integral $ [T(x)]^*(\mu) $ belongs to $Val(X)$.
Since $Val(X)$ is complete (it is a Banach space), the integral is well defined as a limit of Riemann sums. \qed

\item  Given a linear map $F\colon X\to Y$. Now let us define the pull-back map
\begin{eqnarray}\label{E:F-star-sec-2}
F^*\colon Val_k(Y)\to Val_k(X)
\end{eqnarray}
as follows. Let $f\in  Val_k(Y)$. For any subspace $E\in Gr_k^X$ let us define
\begin{eqnarray}\label{E:F-star-sections}
(F^*f)(E)=\left\{\begin{array}{ccc}
                  \ct(F|_E,f(F(E)))&\mbox{ if }& \dim F(E)=k\\
                  0&\mbox{ otherwise }&
                \end{array}\right. .
\end{eqnarray}

\begin{proposition}\label{CL:contin-pull-back}
$F^*f$ is a continuous section of $\cl_k^X$.
\end{proposition}

We will need the following elementary lemma whose proof is easy and is left to the reader. For vector spaces $K,X$ we denote by $Inj(K,X)$ the space of linear imbeddings $K\inj X$.
\begin{lemma}\label{L:contin-lemma-tech}
Let $K$ be a $k$-dimensional vector space. Let $f$ be a not necessarily continuous section of $\cl_k^X$ over $Gr_k^X$. Then $f$ is continuous if and only if the map
$Inj(K,X)\to D(K)$ given by $h\mapsto h^*f$ is continuous.
\end{lemma}

We will need one more lemma.
\begin{lemma}\label{L:contin-l-tech2}
Let $K$ and $X$ be finite dimensional vector spaces over $\FF$, and $\dim K=k$. Then the map
\begin{eqnarray}\label{E:map-leb-meas}
Hom(K,X)\times Val_k(X)\to D(K)
\end{eqnarray}
given by $$(f,\phi)\mapsto f^*\phi\in Val_k(K)=D(K)$$
is jointly continuous.
\end{lemma}
{\bf Proof.} The map (\ref{E:map-leb-meas}) is linear with respect to the second variable, $Val_k(X)$ and $D(K)$ are Banach spaces, and $Hom(K,X)$ is a metrizable locally compact space. Hence it suffices to show that the map
(\ref{E:map-leb-meas}) is separately continuous (see \cite{treves}, Theorem 34.1).

First, let us fix $f\in Hom(K,X)$. The map $Val_k(X)\to D(K)$ given by $\phi\mapsto f^*\phi$ is obviously continuous by the definition of $f^*$.

Next, let us fix $\phi\in Val_k(X)$ and prove that the map $Hom(K,X)\to D(K)$ given by $f\mapsto f^*\phi$ is continuous. It suffices to show that if a sequence of injective maps $\{f_i\}_{i=1}^\infty\subset Hom(K,X)$
converges to a map $f\in Hom(K,X)$, then
$$f_i^*\phi\to f^*\phi.$$
If $f$ is also injective, this follows from Lemma \ref{L:contin-lemma-tech}. Hence, let us assume that $f$ is not injective. We have to show that
$$f^*_i\phi\to 0.$$

Choose a lattice $\Lam_K\subset K$. It suffices to show that
\begin{eqnarray}\label{E:limit-1}
(\phi(f_i(K)))(f_i(\Lam_K))\to 0.
\end{eqnarray}
Fix a lattice $\Lam_X\subset X$. For any $E\in Gr_k^X$ we denote by $vol_E$ the only Lebesgue measure on $E$ such that $vol_E(\Lam_X\cap E)=1$. It suffices to show
$$vol_{f_i(K)}(f_i(\Lam_K))\to 0.$$
The subsets $f_i(\Lam_K)\subset X$converge in the Hausdorff metric to the compact subset $f(\Lam_K)\subset Im(f)$. But $\dim Im(f)<k$. This implies the result. \qed

{\bf Proof of Proposition \ref{CL:contin-pull-back}.} Fix a $k$-dimensional vector space $K$. By Lemma \ref{L:contin-lemma-tech} we have to show that $h\mapsto h^*(F^*f)$ is a continuous map
$Inj(K,X)\to D(K)$. It easily follows from the definition that $h^*(F^*f)=(F\circ h)^*f$. The continuity of $h\mapsto (F\circ h)^*f$ follows from the continuity of $f$ and Lemmas \ref{L:contin-lemma-tech} and \ref{L:contin-l-tech2}. \qed

\item \begin{proposition}\label{pull-back-is-continuous}
The map $F^*$ in  (\ref{E:F-star-sec-2}) is continuous.
\end{proposition}
{\bf Proof.} Let us fix an open bounded subset $\cc\subset X$ which contains the origin. The topology on $C(Gr_k^X,\cl_k^X)$ is given by the norm
$$||g||:=\sup_{E\in Gr_k^X}|\int_{E\cap \cc}g(E)|.$$
Fix an open bounded subset $\ck\subset Y$ such that $F(\cc)\subset \ck$.
Then we have
\begin{eqnarray*}
||F^*f||=\sup_{E\in Gr_k^X,\, \dim F(E)=k}|\int_{F(E\cap \cc)}f(F(E))|\leq \\
\sup_{E\in Gr_k^X,\, \dim F(E)=k}|\int_{F(E)\cap \ck}f(F(E))|\leq ||f||.
\end{eqnarray*}
\qed

\item Let $W\overset{G}{\to}X\overset{F}{\to}Y$ be linear maps. The equality $(F\circ G)^*=G^*\circ F^*$ follows directly from the definition of the pull-back.

\item Let us prove now part (4) of Theorem \ref{T:pull-back-thm}. Namely let us assume that $F\colon X\to Y$ is injective.
Let us show that $F^*(Val^\infty(Y))\subset Val^\infty(X)$. We may and will assume that $X\subset Y$. Let us choose a splitting $Y=X\oplus Z$. It induces
a group imbedding $GL(X)\inj GL(Y)$ given by $g\mapsto (g,Id_Z)$. Clearly any $GL(Y)$-smooth vector is $GL(X)$-smooth. The result follows.

\item \begin{proposition}\label{P:pull-back-double-contin}
The map $Hom(X,Y)\times Val_i(Y)\to Val_i(X)$ given by $(F,f)\mapsto F^*f$ is jointly continuous.
\end{proposition}
{\bf Proof.} It is well known (and easily follows from the Banach-Steinhauss theorem) that any separately continuous map
$$T\times A\to B,$$
where $T$ is a metric space and $A,B$ are Banach spaces, and the operator is linear with respect to the second argument, is jointly continuous.

Hence it suffices to show that our map is separately continuous. By Proposition \ref{pull-back-is-continuous} it remains to show that
if $\xi\in C(Gr_i^Y,\cl_i^Y)$  is fixed then the map $Hom(X,Y)\to C(Gr_i^X,\cl_i^X)$ given by $F\mapsto F^*\xi$ is continuous. Let us prove the continuity of this map at certain fixed $F_0\in Hom(X,Y)$.

Let us fix a lattice $L\subset X$. It suffices to show that for any $\eps>0$ there is a neighborhood of $F_0$ in $Hom(X,Y)$ such that for any $F$ from this neighborhood
\begin{eqnarray}\label{E:local-grassmann-0}
\sup_{E\in Gr_i^X} \big|\int_{E\cap L}(F^*\xi)(E)-\int_{E\cap L}(F^*_0\xi)(E)\big|<\eps.
\end{eqnarray}
It suffices to prove this statement locally, i.e. to prove that any $E_0\in Gr_i^X$ has a neighborhood $\cu_{E_0}\subset Gr_i^X$, and there is a neighborhood $\cv_{E_0}$ of $F_0$ in $Hom(X,Y)$ such that
\begin{eqnarray}\label{E:local-grassmann}
\sup_{E\in \cu_{E_0}} \big|\int_{E\cap L}(F^*\xi)(E)-\int_{E\cap L}(F_0^*\xi)(E)\big|<\eps \mbox{ for any } F\in \cv_{E_0}.
\end{eqnarray}
Indeed then we could choose a finite subcovering $\{\cu_{E_\alp}\}$ of $Gr_i^X$. Then for any  $F\in \cap_\alp \cv_\alp$ one would have (\ref{E:local-grassmann-0}).

Thus let us fix $E_0\in Gr_i^X$.

\hfill

\underline{Case 1.} Let us assume that $\dim F_0(E_0)=i$. By Theorem \ref{T:free-lattice}(i) there exists a linear isomorphism
$H_0\colon \FF^i\tilde\to E_0$ such that $H_0(\co^i)=E_0\cap L$. It suffices to show that for any $\eps>0$ there exists a neighborhood $\cx$ of $H_0$ in $Hom(\FF^i,X)$ and a neighborhood
$\cv$ of $F_0$ in $Hom(X,Y)$ such that for any $H\in \cx$ and any $F\in \cv$ one has
\begin{eqnarray*}
\big|\int_{H^{-1}(Im(H)\cap L)}(F\circ H)^*\xi-\int_{H^{-1}(Im(H)\cap L)}(F_0\circ H)^*\xi\big|<\eps.
\end{eqnarray*}
Note that $H_0^{-1}(Im(H_0)\cap L)=\co^i$. Then we can choose the neighborhood $\cx$ of $H_0$ as in Proposition \ref{P:close-norms}, i.e. so that
$$H^{-1}(Im(H)\cap L)=\co^i \mbox{ for all } H\in \cx.$$
Hence one has to show that
\begin{eqnarray*}
\big|\int_{\co^i}(F\circ H)^*\xi-\int_{\co^i}(F_0\circ H)^*\xi\big|<\eps.
\end{eqnarray*}
But this is clear from the definition of topology on $\cl_i^Y$.

\hfill

\underline{Case 2.} Let us assume that $\dim F_0(E_0)<i$. By Theorem \ref{T:free-lattice} and Remark \ref{Rem:basis-of-lattice} there is an isomorphism
$H_0\colon \FF^i\tilde\to E_0$ such that $H_0(\co^i)=E_0\cap L$ and
$$H_0(\co^k\times \{0_{i-k}\})=(KerF_0)\cap E_0\cap L,$$
where $k:=\dim ((Ker F_0)\cap E_0).$ It suffices to show that for any $\eps>0$ there exists a neighborhood $\cx$ of $H_0$ in $Hom(\FF^i,X)$ and a neighborhood
$\cv$ of $F_0$ in $Hom(X,Y)$ such that for any $H\in \cx$ and any $F\in \cv$ one has
\begin{eqnarray}\label{E:ooo1}
\big|\int_{H^{-1}(Im(H)\cap L)}(F\circ H)^*\xi\big|<\eps.
\end{eqnarray}

Note that $H_0^{-1}(Im(H_0)\cap L)=\co^i$. By Proposition \ref{P:close-norms} $H_0$ has a neighborhood $\cx\subset Hom(\FF^i,X)$ so that
\begin{eqnarray*}
H^{-1}(Im(H)\cap L)=\co^i \mbox{ for all } H\in \cx,
\end{eqnarray*}
Hence (\ref{E:ooo1}) is equivalent to saying that there exists a neighbourhood $\cx$ of $H_0$ and $\cv$ of $F_0$ such that
\begin{eqnarray}\label{E:ooo2}
\big|\int_{\co^i}(F\circ H)^*\xi\big|<\eps \mbox{ for any } H\in \cx, \,F\in \cv.
\end{eqnarray}
It suffices to show that given the linear map $g_0\colon \FF^i\to Y$ such that $rk(g_0)<i$ then there is a neighborhood $\cu\subset Hom(\FF^i,Y)$ of $g_0$ such that
$$\big|\int_{\co^i}g^*\xi\big|<\eps \,\,\, \forall g\in \cu.$$
We leave to the reader this simple and elementary estimate. \qed


\end{paragraphlist}

\section{Construction of Fourier transform on valuations.}\label{S:fourier-transform}
\begin{paragraphlist}
\item Let $\FF$ be a non-Archimedean local field. In this section we construct a $GL(V)$-equivariant isomorphism between spaces of valuations on a vector space $V$ and its dual.

We denote by $V^\vee$ (rather than $V^*$) the dual space of $V$. Nevertheless we will keep $\ast$ to denote duals of one dimensional spaces (thus $D(V)^*$ denotes the dual space of the space $D(V)$ of
$\CC$-valued Lebesgue measures on $V$).

\hfill For a vector subspace $E\subset V$ let us denote by $E^\perp \subset V^\vee$ its annihilator defined by
$$E^\perp:=\{f\in V^\vee|\,\, f(E)=0\}.$$
This induces a $GL(V)$-equivariant homeomorphism
$$Gr_i^V\overset{\perp}{\tilde\to} Gr_{n-i}^{V^\vee}.$$

\item Recall that the fiber of the line bundle $\cl_i^V$ over $E\in Gr_i^V$ is equal, by the definition, to the space $D(E)$ of complex valued Lebesgue  measures on $E$.
We are going to construct a $GL(V)$-equivariant homeomorphism $a_i\colon\cl_i^V\to \cl_{n-i}^{V^\vee}\otimes D(V^\vee)^*$ such that the diagram

\begin{tikzcd}\label{diag-fourier}
\cl_i^V\arrow{r}{a_i}\arrow{d}&\cl_{n-i}^{V^\vee}\otimes D(V^\vee)^*\arrow{d}\\
Gr_i^V\arrow{r}{\perp}&Gr_{n-i}^{V^\vee}
\end{tikzcd}

is commutative and the map $a_i$ is linear on fibers of the bundles (the vertical arrows are the obvious bundle projections).

Fix a linear subspace $E\subset V$. Note that canonically the dual space $E^\vee=V^\vee/E^\perp$. We have the isomorphisms from Section \ref{S:measures}:
\begin{eqnarray}\label{E:isom01}
D(E)\tilde\to D(E^\vee)^* \tilde\to D(V^\vee/E^\perp)^* \tilde\to (D(V^\vee)\otimes D(E^\perp)^*)^* \tilde\to D(E^\perp)\otimes D(V^\vee)^*.
\end{eqnarray}
This defines the required map $a_i$ on the fiber $D(E)$ over $E$.

\item\label{item-fourier-isomorphism} By taking sections of the bundles in the diagram from paragraph \ref{diag-fourier} we get a $GL(V)$-equivariant isomorphism of Banach spaces
\begin{eqnarray}\label{E:fourier-large}
\FF\colon Val_i(V)\tilde\to Val_{n-i}(V^\vee)\otimes D(V^\vee)^*.
\end{eqnarray}
\begin{remark}
For $i=0$ the Fourier transform is the obvious isomorphism $$\FF\colon \CC\to D(V^\vee)\otimes D(V^\vee)^* .$$
For $i=n=\dim V$  the Fourier transform is the obvious isomorphism
$$\FF\colon D(V)\to \underset{\CC}{\underbrace{Val_0(V^\vee)}}\otimes D(V^\vee)^* .$$
\end{remark}

\item We have the following Plancherel type inversion formula.
\begin{theorem}\label{T:plancherel}
The composition
\begin{eqnarray*}
Val_i(V)\overset{\FF_V}{\to} Val_{n-i}(V^\vee)\otimes D(V^\vee)^*\overset{\FF_{V^\vee}\otimes Id}{\to}\\
\left(Val_i(V)\otimes D(V)^*\right)\otimes D(V^\vee)^*=
Val_i(V)
\end{eqnarray*}
is the identity map.
\end{theorem}
{\bf Proof.} First, taking the annihilator $\perp$ twice is the identity map on $Gr_i^V$. Next, for $E\in Gr_i^V$ the composition of the following isomorphisms is the identity map of $D(E)$
\begin{eqnarray*}
D(E)\tilde\to D(E^\vee)^*\tilde\to (D(V^\vee/E^\perp))^*\tilde\to
D(E^\perp)\otimes D(V^\vee)^*\tilde\to\\ D((E^\perp)^\vee)^*\otimes D(V^\vee)^*\tilde\to D(V/E)^*\otimes D(V^\vee)^*\tilde\to
\left(D(E)\otimes D(V)^*\right)\otimes D(V^\vee)^*\tilde\to\\
 D(E).
\end{eqnarray*}
Hence $a_{n-i}\circ a_i=Id$.

\item The two simplest examples of computation of the Fourier transform of valuations are as follows. Let $\mu\in D(X)=Val_{\dim X}(X)$. Then
\begin{eqnarray}\label{E:fourier-density}
\FF(\mu)=\chi_{X^\vee}\otimes \mu.
\end{eqnarray}
The Fourier transform of the Euler characteristic $\chi_X\in Val_0(X)$ can also be easily computed. Let us fix a non-vanishing Lebesgue measure  $vol_X\in D(X)$.
Let $vol_X^{-1}\in D(X)^*\simeq D(X^\vee)$
be the inverse Lebesgue measure. Then
\begin{eqnarray}\label{E:fourier-euler}
\FF(\chi_X)=vol_X^{-1}\otimes vol_X.
\end{eqnarray}

\item Let us give one more example of computation of the Fourier transform.
\begin{lemma}\label{L:Four-onto-comp}
Let $f\colon X\twoheadrightarrow Y$ be a linear map onto. Let $\mu\in D(Y)$. Denote by $q\colon X^\vee \twoheadrightarrow X^\vee/Y^\vee$ the natural quotient map. Then the Fourier transform $\FF(f^*\mu)$ is equal to
$$q^*(vol_{X^\vee/Y^\vee})\otimes (\mu\otimes vol_{X^\vee/Y^\vee}^{-1}),$$
where $vol_{X^\vee/Y^\vee}$ is any non-zero Lebesgue measure on $X^\vee/Y^\vee$. Here $$\mu\otimes vol_{X^\vee/Y^\vee}^{-1}\in D(X).$$
\end{lemma}
A proof is left to the reader.

\end{paragraphlist}

\section{Push-forward on valuations.}\label{S:push-forward}
\begin{paragraphlist}
\item Let $F\colon X\to Y$ be a linear map. We will define a continuous linear map
\begin{eqnarray}\label{E:ppp1}
F_*\colon Val(X)\otimes D(X)^*\to Val(Y)\otimes D(Y)^*
\end{eqnarray}
called the push-forward map and describe it more explicitly in the two cases of injective and surjective maps.

Let us define first the push-forward
$$F_*\colon C(Gr_k^X,\cl_k^X\otimes D(X)^*)\to C(Gr_{k-\dim X+\dim Y}^Y,\cl_{k-\dim X+\dim Y}^Y\otimes D(Y)^*)$$
as follows. We have the dual map $F^\vee\colon Y^\vee\to X^\vee$. Define
\begin{eqnarray}\label{E:def-push-forward}
F_*:=\FF_{Y^\vee}\circ (F^\vee)^*\circ \FF^{-1}_{X^\vee},
\end{eqnarray}
where $\FF_{X^\vee},\FF_{Y^\vee}$ are the Fourier transforms on $X^\vee,Y^\vee$ respectively.
\begin{lemma}
$F_*$ is a continuous map $$Val_k(X)\otimes D(X)^*\to C(Gr_{k-\dim X+\dim Y}^Y,\cl_{k-\dim X+\dim Y}^Y\otimes D(Y)^*).$$
\end{lemma}
This easily follows from the properties of pull-back and Fourier transform.

\item Let us describe the push-forward map when $F\colon X\to Y$ is an imbedding. The description is contained in Propositions \ref{E:ppp5} and \ref{E:ppp6} below.

We will identify $X$ with its image in $Y$.  Denote $c:=\dim Y-\dim X$ the codimension of $X$.

Let $\xi\in Val_k(X)\otimes D(X)^*$. Let $E\in Gr_{k+c}^Y$. We are going to describe $(F_*\xi)(E)$. We will consider two cases: $Y\ne E+X$ and $Y= E+X$.

\begin{proposition}\label{E:ppp5}
 Assume $Y\ne E+X$.  Then  $(F_*\xi)(E)=0$.
\end{proposition}
{\bf Proof.} Let us introduce a notation. For a linear subspace $L\subset V$  consider the isomorphism
\begin{eqnarray}\label{E:ppp5.3}
\alp_{L,V}\colon D(L)\tilde\to D(L^\perp)\otimes D(V^\vee)^*
\end{eqnarray}
which is  the composition of natural isomorphisms (\ref{E:isom01}) with $E$ replaced with $L$:
\begin{eqnarray*}\label{E:ppp5.5}
D(L)\tilde\to D(L^\vee)^* \tilde\to D(V^\vee/L^\perp)^* \tilde\to (D(V^\vee)\otimes D(L^\perp)^*)^* \tilde\to D(L^\perp)\otimes D(V^\vee)^*.
\end{eqnarray*}
We have
\begin{eqnarray}\label{E:ppp7}
(F_*\xi)(E)=\left((\FF_{Y^\vee}\circ F^{\vee \ast}\circ \FF^{-1}_{X^\vee})(\xi)\right)(E)=\alp_{E^\perp, Y^\vee}\left(((F^{\vee\ast}\circ\FF^{-1}_{X^\vee})(\xi))(E^\perp)\right).
\end{eqnarray}
We have to show that the last expression vanishes. It suffices to show that $Ker(F^\vee\colon E^\perp\to X^\vee)\ne 0$. By duality this is equivalent to
$$CoKer(X\to (E^\perp)^\vee)\ne 0.$$
But $(E^\perp)^\vee=Y/(E^\perp)^\perp=Y/E$. Hence, equivalently we have
$$CoKer(X\to Y/E)\ne 0,$$
where the map is the composition of the natural maps $X\to Y\to Y/E$. This is equivalent to our assumption $Y\ne E+X$. \qed

\item Assume now that $Y=E+X$. Denote $E_0=E\cap X$. Then $\dim E_0=k$. The natural map $E\to Y/X$ induces the isomorphism
\begin{eqnarray}\label{E:ppp2.5}
E/E_0\tilde\to Y/X.
\end{eqnarray}

\begin{proposition}\label{E:ppp6}
Assume that $Y=E+X$. Then $(F_*\xi)(E)$ equals the image of $\xi(E_0)\in D(E_0)\otimes D(X)^*$ under the composition of natural isomorphisms from Section \ref{S:measures}
\begin{eqnarray*}\label{E:ppp3}
 D(E_0)\otimes D(X)^*\tilde\to  D(E_0)\otimes D(Y/X)\otimes D(Y)^*\tilde\to\\
  D(E_0)\otimes D(E/E_0) \otimes D(Y)^*\tilde\to
 D(E)\otimes D(Y)^*,
\end{eqnarray*}
where in the second isomorphism we used the isomorphism (\ref{E:ppp2.5}).
\end{proposition}
{\bf Proof.} We can choose a decomposition
\begin{eqnarray}\label{E:decomp-00}
Y=Z\oplus E_0\oplus E_1
\end{eqnarray}
such that
\begin{eqnarray*}
E=E_0\oplus E_1,\,\,
X=Z\oplus E_0.
\end{eqnarray*}

By (\ref{E:ppp7}) we have
\begin{eqnarray}\label{E:ppp8}
(F_*\xi)(E)=\alp_{E^\perp, Y^\vee}\left(((F^{\vee\ast}\circ\FF^{-1}_{X^\vee})(\xi))(E^\perp)\right).
\end{eqnarray}
Recall that $\alp_{E^\perp, Y^\vee}\colon D(E^\perp)\tilde\to D(E)\otimes D(Y)^*$. Using decomposition (\ref{E:decomp-00}) we clearly have $E^\perp =Z^\vee$.  Under this identification we have
$$\alp_{E^\perp, Y^\vee}\colon D(Z^\vee)\to D(E)\otimes D(Y)^*.$$
The dual map $F^\vee$ is the natural projection $Z^\vee\oplus E_0^\vee\oplus E_1^\vee\to Z^\vee\oplus E_0^\vee$.
The subspace $E^\perp=Z^\vee$ is mapped identically under these identifications. Hence it follows
\begin{eqnarray*}
((F^{\vee\ast}\circ\FF^{-1}_{X^\vee})(\xi))(E^\perp)=\left(F^{\vee\ast}(\FF^{-1}_{X^\vee}\xi)\right)(Z^\vee)=(\FF^{-1}_{X^\vee}\xi)(Z^\vee)=\alp^{-1}_{Z^\vee,X^\vee} (\xi(E_0)),
\end{eqnarray*}
where $\alp_{Z^\vee,X^\vee}\colon D(Z^\vee)\to D(E_0)\otimes D(X)^*.$
Substituting this into (\ref{E:ppp8}) we get
\begin{eqnarray*}
(F_*\xi)(E)=\left(\alp_{E^\perp, Y^\vee}\circ \alp^{-1}_{Z^\vee,X^\vee}\right) (\xi(E_0)).
\end{eqnarray*}
The map $\alp_{E^\perp, Y^\vee}\circ \alp^{-1}_{Z^\vee,X^\vee}\colon D(E_0)\otimes D(X)^*\to D(E)\otimes D(Y)^*$ coincides with the claimed one. \qed

\item \begin{proposition}\label{P:push-forward-onto}
Let us assume that a linear map $F\colon X\to Y$ is onto. Then push-forward of a smooth valuation is smooth:
$$F_*(Val^\infty(X)\otimes D(X)^*)\subset Val^\infty(Y)\otimes D(Y)^*.$$
\end{proposition}
{\bf Proof.} Recall that $F_*=\FF_{Y^\vee}\circ (F^\vee)^*\circ \FF^{-1}_{X^\vee}$. Since $F^\vee$ is an imbedding, $(F^\vee)^*$ maps
smooth valuations to smooth ones by Theorem \ref{T:pull-back-thm}(4). By Section \ref{S:fourier-transform}, paragraph \ref{item-fourier-isomorphism}, $\FF$ is an isomorphism
between spaces of smooth valuations. The result follows. \qed

\item In this paragraph we will describe explicitly the push-forward map when $F\colon X\to Y$ is onto. Denote $K:=Ker(F)$. Clearly $X/K\simeq Y$.
\begin{proposition}\label{claim-ppp9}
Let $F\colon X\to Y$ be onto. Let $\xi\in Val_k(X)\otimes D(X)^*$. Let $E\in Gr_{k-\dim X+\dim Y}^Y$ with $k\geq \dim X-\dim Y$. Denote $\tilde E:=F^{-1}(E)$.
Then $(F_*\xi)(E)$ is equal to the image of $\xi(\tilde E)$ under the composition of natural isomorphisms
\begin{eqnarray*}
D(\tilde E)\otimes D(X)^*\tilde \to D(E)\otimes D(K)\otimes D(X)^*\\
\tilde\to D(E)\otimes D(K)\otimes D(K)^*\otimes D(Y)^*\tilde\to D(E)\otimes D(Y)^*,
\end{eqnarray*}
where we used isomorphism $\tilde E/K\simeq E$.
\end{proposition}
{\bf Proof.} Let us choose a splitting $X=K\oplus Y$. Then $F^\vee\colon Y^\vee\to K^\vee\oplus Y^\vee$ is the obvious imbedding $y\mapsto (0,y)$. Also $\tilde E=K\oplus E$.
We will identify the annihilator $E^\perp$ of $E$ in $Y^\vee$ with its image in $X^\vee$ and denote in the same way by the abuse of notation.

We have
\begin{eqnarray*}
(F_*\xi)(E)=\left((\FF_{Y^\vee}\circ F^{\vee \ast}\circ \FF_{X^\vee}^{-1})(\xi)\right)(E)=\alp_{E^\perp,Y^\vee}\left(((F^{\vee\ast}\circ \FF_{X^\vee}^{-1})(\xi))(E^\perp)\right)=\\
\alp_{E^\perp,Y^\vee}\left((\FF^{-1}_{X^\vee}\xi)(E^\perp)\right)=\alp_{E^\perp,Y^\vee}\left(\alp_{E^\perp,X^\vee}^{-1}(\xi(\tilde E))\right)=\\
(\alp_{E^\perp,Y^\vee}\circ \alp_{E^\perp,X^\vee}^{-1})(\xi(\tilde E)),
\end{eqnarray*}
where we recall that
\begin{eqnarray*}
\alp_{E^\perp,X^\vee}\colon D(E^\perp)\tilde\to D(\tilde E)\otimes D(X)^*,\\
\alp_{E^\perp,Y^\vee}\colon D(E^\perp)\tilde\to D(E)\otimes D(Y)^*.
\end{eqnarray*}
It is easy to see that the map $\alp_{E^\perp,Y^\vee}\circ \alp_{E^\perp,X^\vee}^{-1}$ coincides with the map from the proposition. \qed

\end{paragraphlist}

\section{Exterior product on valuations.}\label{S:Ext-product}
\begin{paragraphlist}
\item\label{item:ext-prod-mention} The goal of this section is to construct a bilinear map for any two finite dimensional vector spaces $X$ and $Y$ over a non-Archimedean local field $\FF$
\begin{eqnarray}\label{E:ooo3}
\boxtimes\colon Val_i(X)\times Val^\infty_j(Y)\to Val_{i+j}(X\times Y)
\end{eqnarray}
which is $GL(X)\times GL(Y)$-equivariant and continuous with respect to the first variable when the second one is fixed.

\item\label{item:exter-densities} First let us construct the map (\ref{E:ooo3}) in the special case $Val^\infty_j(Y)=D(Y)$. In this case we need to construct a continuous bilinear map
\begin{eqnarray}\label{E:ooo4}
\boxtimes \colon Val_i(X)\otimes D(Y)\to Val_{i+\dim Y}(X\times Y).
\end{eqnarray}
This map is essentially the push-forward map for the obvious imbedding $$F\colon X\inj X\times Y$$
given by $F(x)=(x,0)$.
Indeed, the push-forward map for this imbedding, as defined in Section \ref{S:push-forward}, is a continuous linear map
$$F_*\colon Val_i(X)\otimes D(X)^*\to  Val_{i+\dim Y}(X\times Y)\otimes D(X\times Y)^*.$$
Since $D(X\times Y)^*=D(X)^*\otimes D(Y)^*$ we get the required map (\ref{E:ooo4}) by twisting all spaces by $D(X\times Y)$.

Exterior product of $\phi\in C(Gr_i^X,\cl_i), \mu\in D(Y)$ is denoted by $\phi\boxtimes \mu$. Thus
\begin{eqnarray}\label{E:yyyyy}
\phi\boxtimes \mu=(F_*\otimes Id_{D(X\times Y)})(\phi\otimes\mu).
\end{eqnarray}

\item\label{item:ext-dens-val}
\begin{lemma}\label{L:ext-dens-pull-back-dens}
Let $T\colon X\to Z$ be a linear map. Let $\mu\in D(Z),\, \nu\in D(W)$. Then
$$(T^*\mu)\boxtimes \nu=(T\times Id_W)^*(\mu\boxtimes \nu).$$
\end{lemma}
{\bf Proof.} If $T$ is not onto then both sides vanish. Thus let us assume that $T$ is onto. Fix a linear subspace $E\subset X\times W$ with $\dim E=\dim Z+\dim W$.
We have to show that
\begin{eqnarray*}
\left((T^*\mu)\boxtimes \nu\right)(E)=\left((T\times Id)^*(\mu\boxtimes \nu)\right)(E).
\end{eqnarray*}
Both sides are continuous in $E$. Hence it suffices to assume that $E$ is generic.
Denote $E_0:=E\cap X$. For generic $E$ the maps $(T\times Id_W)|_E\colon E\to Z\times W$  and $T_{E_0}\colon E_0\to Z$ are isomorphisms. The lemma follows from the commutativity of the following diagram
\newline
\begin{tikzcd}
D(Z\times W)\simeq  D(Z)\otimes D(W)\arrow[r,"(T\times Id_W)^*"]\arrow[d,"T^*\otimes Id_W"]&D(E)\arrow[d]\\
D(E_0)\otimes D(W)\arrow[r]&D(E)\otimes D(E/E_0)^*\otimes D(W)
\end{tikzcd}.

\qed


\item Let us denote $n=\dim Y$. To construct the map (\ref{E:ooo3}) in general, first we will construct, using the construction from paragraph \ref{item:exter-densities}, a bilinear map
\begin{eqnarray}\label{E:ooo6}
\tilde \boxtimes \colon Val_i(X)\times C^\infty(Gr_{n-j}^Y,\cm_{n-j}^Y)\to Val_{i+j}(X\times Y),
\end{eqnarray}
which is continuous with respect to the first variable,  where the bundle $\cm_{n-j}^Y$ was defined in Section \ref{S:space of valuations}, paragraph \ref{item:m-l}.

Let $\xi$ be a section of $\cm_{n-j}^Y$.
Recall that for any $F\in Gr_{n-j}^Y$ one has $\xi(F)\in D(Y/F)\otimes |\ome_{n-j}^Y|\big|_F$, where $|\ome_{n-j}^Y|$ is the line bundle of densities over $Gr_{n-j}^Y$. Then
by the construction of paragraph \ref{item:exter-densities}
$$\phi\boxtimes \xi(F)\in Val_{i+j}(X\times (Y/F))\otimes |\ome_{n-j}^Y|\big|_F.$$
Let us define
\begin{eqnarray}\label{E:ext-prod-gen99}
\phi\tilde\boxtimes \xi:=\int_{F\in Gr_{n-j}^Y}(Id_X\times p_F)^*(\phi\boxtimes\xi(F)) \in Val_{i+j}(X\times Y),
\end{eqnarray}
where $p_F\colon Y\to Y/F$ is the canonical projection. Note that for the expression under the integral one has
$$(Id_X\times p_F)^*(\phi\boxtimes\xi(F))\in Val_{i+j}(X\times Y)\otimes |\ome_{n-j}^Y|\big|_F.$$

We will prove the following
\begin{claim}\label{CL:well-defined-exter-product}
The integral in (\ref{E:ext-prod-gen99}) is
\newline
(1) well defined;
\newline
(2) linear with respect to $\phi$ and $\xi$;
\newline
(3) continuous with respect to $\phi$.
\end{claim}

\item {\bf Proof of Claim \ref{CL:well-defined-exter-product}.} Let us fix a smooth positive measure $\nu_0$ on $Gr_{n-j}^Y$. Then we can write uniquely $\xi=\tilde\xi\otimes \nu_0$ where $\tilde\xi(F)\in D(Y/F)$, and $\tilde\xi$ is a continuous section.
Then
$$(Id_X\times p_F)^*(\phi\boxtimes \xi(F))=(Id_X\times p_F)^*(\phi \boxtimes \tilde\xi(F))\cdot\nu_0,$$
and $\phi\tilde\boxtimes \xi=\int_{Gr_{n-j}^Y}(Id_X\times p_F)^*(\phi\boxtimes\tilde\xi(F))d\nu_0(F)$, where $$(Id_X\times p_F)^*(\phi\boxtimes\tilde\xi(F))\in Val_{i+j}(X\times Y).$$
It suffices to show that the map $Val_i(X)\times Gr_{n-j}^Y\to Val_{i+j}(X\times Y)$  given by
\begin{eqnarray}\label{E:lll-1}
(\phi,F)\mapsto (Id_X\times p_F)^*(\phi\boxtimes\tilde\xi(F))
\end{eqnarray}
is continuous.
For let us fix $F_0\in Gr_{n-j}^Y$. It has a neighborhood $\cu\subset Gr_{n-j}^Y$ over which there is a trivialization of the vector bundle
whose fiber over $F$ is $Y/F$, thus this bundle is isomorphic to $\cu\times \FF^j$. Under this identification $\tilde\xi \colon \cu\to D(\FF^j)$ is a continuous map. Hence the map $(\phi,F)\mapsto \phi\boxtimes\tilde\xi(F)$ is a continuous map
$Val_i(X)\times\cu\to Val_{i+j}(X\times \FF^j)$ by paragraph \ref{item:exter-densities}.

Under the above trivialization the map $p_F$ becomes a linear map $Y\to \FF^j$ which we denote in the same way; it depends continuously on $F$. By
Proposition \ref{P:pull-back-double-contin} the expression under the integral (\ref{E:ext-prod-gen99}) is jointly continuous with respect to $(\phi,F)$. Hence parts (1), (3) follow. Part (2) is obvious.  \qed

\item Recall that in Section \ref{S:space of valuations} we defined a linear isomorphism
$$\cd\colon C^\infty(Gr^Y_{n-j},\cm^Y_{n-j})\tilde\to C^\infty(Gr^Y_j,\cl_j^Y)=Val_j^\infty(Y).$$
Define now the exterior product
$$\boxtimes\colon Val_i(X)\times Val_j^\infty(Y)\to Val_{i+j}(X\times Y)$$
to be
\begin{eqnarray}\label{Def:exter-product}
\boxtimes=\tilde\boxtimes \circ(Id_{Val_i(X)}\times \cd^{-1}).
\end{eqnarray}

Claim \ref{CL:well-defined-exter-product} immediately implies that $\boxtimes$ is a bilinear map continuous with respect to the first argument when the second one is fixed.

\item Let us generalize Lemma \ref{L:ext-dens-pull-back-dens} as follows.
\begin{proposition}\label{P:pull-back-exterior-product}
Let $T\colon X\to Z$ be a linear map. Let $W$ be another vector space. Let $\ome\in Val_i(Z),\, \xi\in Val_j^\infty(W)$.  Then
$$T^*\ome\boxtimes \xi=(T\times Id_W)^*(\ome\boxtimes \xi).$$
\end{proposition}
{\bf Proof.} Since the pull-back $T^*\colon Val_i(Z)\to Val_i(X)$
is continuous by Theorem \ref{T:pull-back-thm} and the exterior product is continuous with respect to the first argument by Claim \ref{CL:well-defined-exter-product}(3), Lemma \ref{L:pull-back-density}
implies that we may assume that
$$\ome=S^*\mu,$$
where $S\colon Z\to \FF^i$ is a linear map and $\mu\in D(\FF^i)$. We can represent $\xi$ in the form
$$\xi=\int_{F\in Gr^W_{n-j}}p_F^*(\zeta(F)),$$
where $n=\dim W$, $p_F\colon W\to W/F$ is the canonical quotient map,  $\zeta$ is a smooth section of the line bundle $\cm_{n-j}^W$ over Grassmannian $Gr_{n-j}^W$ whose fiber over $F$ is $D(W/F)$ tensorized with the
fiber of the line bundle of densities on $Gr_{n-j}^W$.

We have
\begin{eqnarray*}
T^*\ome\boxtimes \xi\overset{(\ref{E:ext-prod-gen99})}{=}\\
\int_{F\in Gr}(Id\times p_F)^*(T^*\ome\boxtimes \zeta(F))=
\int_{F\in Gr} (Id\times p_F)^*(T^*S^*\mu\boxtimes \zeta(F))\overset{Lemma\, \ref{L:ext-dens-pull-back-dens}}{=}\\
\int_{F\in Gr}(Id\times p_F)^*(ST\times Id)^*(\mu\boxtimes \zeta(F))\overset{Thm \,\ref{T:pull-back-thm}(2)}{=}\\
\int_{F\in Gr}(ST\times p_F)^*(\mu\boxtimes \zeta(F))\overset{Thm\,\ref{T:pull-back-thm}}{=}\\
\int_{F\in Gr} (T\times Id)^*(Id\times p_F)^*(S\times Id)^*(\mu\boxtimes \zeta(F))\overset{Lemma\, \ref{L:ext-dens-pull-back-dens}}{=}\\
(T\times Id)^*\int_{F\in Gr}(Id\times p_F)^*(\ome\boxtimes \zeta(F))\overset{(\ref{E:ext-prod-gen99})}{=}\\
(T\times Id)^*(\ome\boxtimes \xi).
\end{eqnarray*}
\qed


\item \begin{lemma}\label{L:gen-associativity}
Let $X,Y,Z$ be finite dimensional $\FF$-vector spaces. Then the two maps
$$Val^\infty(X)\times Val(Y)\times Val^\infty(Z)\to Val(X\times Y\times Z)$$
given respectively by
\begin{eqnarray*}
(\phi,\psi,\xi)\mapsto (\phi\boxtimes \psi)\boxtimes \xi,\\
(\phi,\psi,\xi)\mapsto \phi\boxtimes (\psi\boxtimes \xi)
\end{eqnarray*}
coincide with each other.
\end{lemma}
{\bf Proof.} \underline{Step 1.} Both maps are 3-linear and continuous with respect to $\psi$. By Lemma \ref{L:pull-back-density} linear combinations of valuations of the form $T^*\mu$, where $\mu$ is a Lebesgue measure,
are dense in $Val(Y)$.  Thus it suffices to assume that $\psi=T^*\mu$ where $T\colon Y\to W$, $\mu \in D(W)$.

\underline{Step 2.} Let us assume that $\psi=T^*\mu$. Applying twice Proposition \ref{P:pull-back-exterior-product} we have
\begin{eqnarray*}
(\phi\boxtimes T^*\mu)\boxtimes \xi=\left((Id\times T)^*(\phi\boxtimes\mu)\right)\boxtimes \xi=\\
(Id\times T\times Id)^*((\phi\boxtimes \mu)\boxtimes \xi).
\end{eqnarray*}
Similarly
$$\phi\boxtimes (T^*\mu\boxtimes \xi)=(Id\times T\times Id)^*(\phi\boxtimes (\mu\boxtimes \xi)).$$
Hence it remains to show that $(\phi\boxtimes \mu)\boxtimes \xi=\phi\boxtimes (\mu\boxtimes \xi)$.

\underline{Step 3.} Thus let us assume that $\psi=\mu$ is a density, in particular is a smooth valuation.
Let us denote $x:=\dim X,\, z:=\dim Z$. Since $\phi$ and $\xi$ are smooth valuations they can be presented as
\begin{eqnarray*}
\phi=\int_{Gr^X_{x-i}}p_F^*(\kappa(F)), \, \xi=\int_{Gr^Z_{z-j}}q_E^*\nu(E),
\end{eqnarray*}
where $p_F\colon X\to X/F,\, \, q_E\colon Z\to Z/E$ are the natural quotient maps, and $\kappa$ and $\nu$ are smooth sections of the line bundles $\cm^X_{x-i}$ over $Gr^X_{x-i}$ and $\cm^Z_{z-i}$ over $Gr^Z_{z-j}$
respectively. Then using Proposition \ref{P:pull-back-exterior-product} several times we get
\begin{eqnarray*}
\phi\boxtimes(\mu\boxtimes \xi)=\left(\int_{F\in Gr^X_{x-i}} p_F^*(\kappa(F))\right)\boxtimes \left(\mu\boxtimes \int_{E\in Gr^Z_{z-j}} q_E^*(\nu(E))\right)=\\
\left(\int_{F\in Gr^X_{x-i}} p_F^*(\kappa(F))\right)\boxtimes\left( \int_{E\in Gr^Z_{z-j}}\left(\mu\boxtimes  q_E^*(\nu(E))\right)\right)\overset{Prop. \,\,\, \ref{P:pull-back-exterior-product}}{=}\\
\left(\int_{F\in Gr^X_{x-i}} p_F^*(\kappa(F))\right)\boxtimes\left( \int_{E\in Gr^Z_{z-j}}(Id_Y\times q_E)^*\left(\mu\boxtimes  \nu(E)\right)\right)\overset{Prop. \,\,\, \ref{P:pull-back-exterior-product}}{=}\\
\int_{E\in Gr^Z_{z-j}}(Id_X\times (Id_Y\times q_E))^* \left(\int_{F\in Gr^X_{x-i}} p_F^*(\kappa(F))\boxtimes\left( \mu\boxtimes  \nu(E)\right)\right)\overset{Prop. \,\,\, \ref{P:pull-back-exterior-product}}{=}\\
\int_{E\in Gr^Z_{z-j}}\int_{F\in Gr^X_{x-i}}(Id_X\times Id_Y\times q_E)^*\circ (p_F\times (Id_Y\times Id_{Z/E}))^*\left(\kappa(F)\boxtimes (\mu\boxtimes \nu(E))\right)\overset{Thm\,\,\, \ref{T:pull-back-thm}}{=}\\
\int_{E\in Gr^Z_{z-j}}\int_{F\in Gr^X_{x-i}}(p_F\times Id_Y\times q_E)^*\left(\kappa(F)\boxtimes (\mu\boxtimes \nu(E))\right).
\end{eqnarray*}
Since the exterior product on Lebesgue measures is associative we have
\begin{eqnarray*}
\phi\boxtimes(\mu\boxtimes \xi)=\\
\int_{E\in Gr^Z_{z-j}}\int_{F\in Gr^X_{x-i}}(p_F\times Id_Y\times q_E)^*\left(\left(\kappa(F)\boxtimes \mu\right)\boxtimes \nu(E)\right)\overset{Thm\,\,\, \ref{T:pull-back-thm}}{=}\\
\int_{E\in Gr^Z_{z-j}}\int_{F\in Gr^X_{x-i}}(p_F\times Id_Y\times Id_Z)^*\circ ((Id_{X/F}\times Id_Y)\times q_E)^*\left(\left(\kappa(F)\boxtimes \mu\right)\boxtimes \nu(E)\right)\overset{Prop. \,\,\, \ref{P:pull-back-exterior-product}}{=}\\
\int_{E\in Gr^Z_{z-j}}\int_{F\in Gr^X_{x-i}}(p_F\times Id_Y\times Id_Z)^*\left(\left(\kappa(F)\boxtimes \mu\right)\boxtimes q_E^*\nu(E)\right)=\\
\int_{F\in Gr^X_{x-i}}(p_F\times Id_Y\times Id_Z)^*\left(\left(\kappa(F)\boxtimes \mu\right)\boxtimes\int_{E\in Gr^Z_{z-i}} q_E^*\nu(E)\right)=\\
\int_{F\in Gr^X_{x-i}}((p_F\times Id_Y)\times Id_Z)^*\left(\left(\kappa(F)\boxtimes \mu\right)\boxtimes\xi\right)\overset{Prop. \,\,\, \ref{P:pull-back-exterior-product}}{=}\\
\left(\int_{F\in Gr^X_{x-i}}(p_F\times Id_Y)^*(\kappa(F)\boxtimes \mu)\right)\boxtimes\xi\overset{Prop. \,\,\, \ref{P:pull-back-exterior-product}}{=}\\
\left(\int_{F\in Gr^X_{x-i}} p_F^*\kappa(F)\boxtimes\mu\right)\boxtimes \xi= (\phi\boxtimes\mu)\boxtimes \xi.
\end{eqnarray*}
\qed

\item Let us prove the following identity.
\begin{proposition}\label{P:ext-prod-1}
Let $\chi_X\in Val_0(X)$, $\chi_Y\in Val_0(Y)$ be the unit valuations (Euler characteristics).
Then $\chi_X\boxtimes \chi_Y=\chi_{X\oplus Y}$.
\end{proposition}
{\bf Proof.} Let $p_X\colon X\to \{0\},\, p_Y\colon Y\to \{0\}$ be the obvious maps to the zero-space.
Then formally $\chi_X=p_X^*\chi_0,\, \chi_Y=p_Y^*\chi_0$, where $\chi_0\in Val_0(\{0\})$ is the Euler characteristic on $\{0\}$. Then
\begin{eqnarray*}
\chi_X\boxtimes \chi_Y=p_X^*\chi_0\boxtimes p_Y^*\chi_0\overset{\mbox{Lemma }\ref{L:ext-dens-pull-back-dens}}{=}\\
(Id_X\times p_Y)^*(p_X^*\chi_0\boxtimes \chi_0)\overset{\mbox{Lemma }\ref{L:ext-dens-pull-back-dens}}{=}\\
(Id_X\times p_Y)^*\circ (p_X\times Id_{\{0\}})^*(\chi_0\boxtimes \chi_0)\overset{\mbox{Thm }\ref{T:pull-back-thm}(2)}{=}\\
(p_X\times p_Y)^*(\chi_0\boxtimes \chi_0)=(p_X\times p_Y)^*(\chi_0)=\chi_{X\oplus Y}.
\end{eqnarray*}

\end{paragraphlist}


\section{Product on smooth valuations.}\label{S:product}
Let $\Delta\colon V\to V\times V$ denote the diagonal imbedding, i.e. $\Delta(v)=(v,v)$. Let $\phi\in Val(V),\psi\in Val^\infty(V)$. Define their product $$Val(V)\times Val^\infty(V)\to Val(V)$$ by
$\phi\cdot\psi:=\Delta^*(\phi\boxtimes\psi).$
The product is a bilinear map continuous with respect to the first valuation when the second one is fixed.
\begin{lemma}
Product of smooth valuations is smooth.
\end{lemma}
{\bf Proof.} The product map $Val(V)\otimes Val^\infty(V)\to Val(V)$ is $GL(V)$-equivariant. Hence it maps smooth vectors to smooth ones. \qed

\hfill

Thus we got a product on smooth valuations $Val^\infty(V)\times Val^\infty(V)\to Val^\infty(V)$.
\begin{theorem}\label{T:product-properties}
Equipped with the above product, $Val^\infty(V)$ is a commutative associative algebra with unit $1\in Val_0(V)=\CC$.
It is graded, namely $Val_i^\infty(V)\cdot Val_j^\infty(V)\subset Val_{i+j}^\infty(V)$.
\end{theorem}
{\bf Proof.} Let us prove commutativity. Let us denote by $\sigma\colon V\times V\to V\times V$ the involution $\sigma(x,y)=(y,x)$.
 Then we have
\begin{eqnarray*}
\psi\cdot \phi =\Delta^*(\psi\boxtimes \phi)=\Delta^*(\sigma^*(\phi\boxtimes \psi))=\\
(\Delta^*\circ \sigma^*)(\phi\boxtimes \psi)\overset{Thm\,\,\,\ref{T:pull-back-thm}(2)}{=}
(\sigma\circ \Delta)^*(\phi\boxtimes \psi)=\Delta^*(\phi\boxtimes \psi),
\end{eqnarray*}
where in the last equality we used that $\sigma\circ \Delta=\Delta$.

Let us prove associativity. We have
\begin{eqnarray*}
(\phi\cdot\psi)\cdot \xi=\Delta^*((\phi\cdot\psi)\boxtimes\xi)=\\
\Delta^*(\Delta^*(\phi\boxtimes \psi)\boxtimes \xi)\overset{Prop. \,\,\, \ref{P:pull-back-exterior-product}}{=}\\
\Delta^*(\Delta\times Id_V)^*((\phi\boxtimes \psi)\boxtimes \xi)\overset{Thm \,\,\, \ref{T:pull-back-thm}(2)}=\\
((\Delta\times Id_V)\circ \Delta)^*((\phi\boxtimes \psi)\boxtimes \xi).
\end{eqnarray*}
It is easy to see that $(\Delta\times Id_V)\circ \Delta=\Delta_3$ where $\Delta_3\colon V\to V\times V\times V$ is given by $\Delta_3(v)=(v,v,v)$.
Thus $$(\phi\cdot\psi)\cdot \xi=\Delta_3^*((\phi\boxtimes \psi)\boxtimes \xi).$$
Similarly
$$\phi\cdot(\psi\cdot \xi)=\Delta_3^*(\phi\boxtimes (\psi\boxtimes \xi)).$$
But $(\phi\boxtimes \psi)\boxtimes \xi=\phi\boxtimes (\psi\boxtimes \xi)$ by Lemma \ref{L:gen-associativity}.

The rest of properties are trivial. \qed

\section{Poincar\'e duality for valuations.}\label{S:poincare}
\begin{paragraphlist}
\item The main result of this section is
\begin{theorem}\label{T:poincare}
For any $0\leq i\leq n$ the bilinear form given by the product on valuations
$$Val_i(V)\times Val_{n-i}^\infty(V)\to Val_n(V)$$
is a non-degenerate pairing, i.e. for any non-zero $\phi\in Val_i(V)$ there exists $\psi\in Val_{n-i}^\infty(V)$ such that $\phi\cdot \psi\ne 0$.
\end{theorem}
\item To prove the theorem, let us observe that the subspace
$$\{\phi\in Val_i(V)|\, \, \phi\cdot \psi =0 \,\forall \psi\in Val_{n-i}^\infty(V)\}$$
is a $GL(V)$-invariant closed subspace. By the irreducibility property (Proposition \ref{P:irred}(2)) it suffices to show that it is not equal to $Val_i(V)$.
This follows from the following slightly more general lemma which also will be needed
later on.
\begin{lemma}\label{L:product-positive}
Let $ i,j\geq 0$ be such that $i+j\leq n$. Let $\phi=\cd(f)$, $\psi=\cd(g)$ where $f\in C^\infty(Gr_{n-i}^V,\cm_{n-i}^V)$ and $g\in C^\infty(Gr_{n-j}^V,\cm_{n-j}^V)$
(see Section \ref{S:space of valuations}, paragraph \ref{item-cD} for the definition of the operator $\cd$) be non-negative sections, both not identically 0.
Then $\phi\cdot \psi\ne 0$, and for any $E\in Gr_{i+j}^V$ one has $(\phi\cdot\psi)(E)\geq 0$.
\end{lemma}
{\bf Proof.} For $F\in Gr_{n-i}^V$ we denote by $p_F\colon V\to V/F$ the quotient map, and similarly for $L\in Gr_{n-j}^V$ we denote by $q_L\colon V\to V/L$ the quotient map.  We have
\begin{eqnarray*}
\phi\boxtimes \psi=\int_{F\in Gr_{n-i}^V}p_F^*f(F)\boxtimes \int_{L\in Gr_{n-j}^V}q_L^*g(L)=\\
\int_{F\in Gr_{n-i}^V}\left(p_F^*f(F)\boxtimes \int_{L\in Gr_{n-j}^V} q_L^*g(L)\right)\overset{Prop.\,\,\, \ref{P:pull-back-exterior-product}}{=}\\
\int_{F\in Gr_{n-i}^V}(p_F\times Id_V)^*\left(f(F)\boxtimes \int_{L\in Gr_{n-j}^V}q_L^*g(L)\right)\overset{Prop.\,\,\, \ref{P:pull-back-exterior-product}}{=}\\
\int_{F\in Gr_{n-i}^V}\int_{L\in Gr_{n-j}^V}(p_F\times Id_V)^*\circ (Id_{V/F}\times q_L)^*(f(F)\boxtimes g(L))\overset{Thm \,\,\, \ref{T:pull-back-thm}(2)}{=}\\
\int_{F\in Gr_{n-i}^V}\int_{L\in Gr_{n-j}^V}(p_F\times q_L)^*(f(F)\boxtimes g(L)).
\end{eqnarray*}
Let  $\Delta\colon V\to V\times V$ be the diagonal imbedding, i.e. $\Delta(v)=(v,v)$.
Let $E\in Gr_{i+j}^V$. Since $\phi\cdot\psi=\Delta^*(\phi\boxtimes\psi)$ we have
\begin{eqnarray}\label{E:star01}
(\phi\cdot \psi)(E)=
\int_{F\in Gr_{n-i}^V}\int_{L\in Gr_{n-j}^V}(p_{E,F\times L})^*(f(F)\boxtimes g(L))\in D(E),
\end{eqnarray}
where $p_{E,F\times L}\colon E\to V/F\times V/L$ is the map $v\mapsto (p_F(v),q_L(v))$ with $p_F,q_L$ being the natural quotient maps. The integrand is non-negative.
 There is $E$ such that the integrand is strictly positive for a pair $(F,L)$. Hence the integral is positive. \qed

\item Below we will use the last lemma to show that powers of certain valuation from $Val_1^\infty(V)$ do not vanish.

Let $V$ be a finite dimensional vector space over a non-Archimedean local field $\FF$. Let us fix a lattice $\Lam\subset V$.
Let $GL(\Lam)$ denote the subgroup of $GL(V)$ consisting of such transformations $T$ such that $T(\Lam)=\Lam$. Thus $GL(\Lam)\simeq GL_n(\co)$.

Let  $V_1\in Val_1^\infty (V)$ be a $GL(\Lam)$-invariant non-zero valuation. Such a valuation is unique up to a proportionality. We normalize $V_1$ so that its restriction to any line $l\subset V$ is the only Lebesgue measure
 whose value on the set $l\cap \Lam$ is equal to 1. This characterizes $V_1$ uniquely.
\begin{proposition}\label{P:powers-mean-width}
 For any $1\leq k\leq n$ and any $E\in Gr_k^V$ one has $(V_1^k)(E)>0$. In particular $V_1^k\ne 0$  for any $1\leq k\leq n$.
\end{proposition}
{\bf Proof.} Let us prove it by induction in $k$. For $k=1$ this is clear. Assume $(V_1^{k-1})(L)> 0$ for any $L\in Gr_{k-1}^V$.
This valuation is $GL_n(\co)$-invariant and hence can be  presented as $V_1^{k-1}=\cd(g)$ where $g$ is $GL_n(\co)$-invariant section of $\cm_{n-k+1}^V$.
$g$ must be positive since $(V_1^{k-1})(L)>0$ for any $L\in Gr_{k-1}^V$. Also $V_1=\cd(f)$ where $f$ is $GL_n(\co)$-invariant. Hence $f>0$.
Hence by Lemma \ref{L:product-positive} and $GL_n(\co)$-invariance the valuation $V_1\cdot V_1^{k-1}=V_1^k$ satisfies the conclusion. \qed

\end{paragraphlist}

\section{Integral transforms on Grassmannians.}\label{S:integral-transform}
\begin{paragraphlist}

\item First let us discuss the Radon transform. Let us fix a lattice $\Lam \subset V$. Let $0\leq p<q\leq n-1$. Define the Radon transform
$$R_{pq}\colon C^\infty(Gr_q^V)\to C^\infty(Gr_p^V)$$
as follows
$$(R_{pq}f)(E)=\int_{\substack{ F\supset E\\ \dim F=q}}f(F)d\nu_E(F),$$
where the integration is over the manifold (isomorphic to $Gr_{q-p}^{V/E}$) of $q$-subspaces $F$ containing $E$ with respect to the only probability Haar 
measure $\nu_E$ which is invariant under all $GL(\Lam)$-transformations preserving $E$.
Clearly the operator $R_{pq}$ is $GL(\Lam)$-equivariant.

\begin{theorem}[\cite{petrov-chernov}]\label{T:radon-tr}
If $p+q=n$ then $R_{pq}\colon C^\infty(Gr_q^V)\to C^\infty(Gr_p^V)$ is an isomorphism.
\end{theorem}

We immediately get
\begin{corollary}\label{CoR:K-types}
If $p+q=n$ then the representations of $GL_n(\co)$ in $C^\infty(Gr_p^V)$ and in $C^\infty(Gr_q^V)$ are isomorphic.
\end{corollary}

\item The goal of this paragraph is to explicitly describe the operator $\cd$ from Section \ref{S:space of valuations}, paragraph \ref{item-cD}.
For any linear subspace $M\subset V$ let us denote by $\mu^M$ the only Lebesgue measure on $M$ whose value on $M\cap \Lam$ is equal to 1.
In particular, $\mu^V$ is the Lebesgue measure on $V$ such that $\mu^V(\Lam)=1$.
Let us denote by $\mu_M$ the only Lebesgue measure on $V/M$ whose value on $\Lam/(M\cap \Lam)\subset V/M$ is equal to 1.

This induces a $GL(\Lam)$-invariant trivialization of the bundle $\cl_i^V$. The $GL(\Lam)$-invariant Haar probability measure on Grassmannians induces $GL(\Lam)$-invariant trivialization
of the line bundle of densities. Hence we get $GL(\Lam)$-invariant trivialization of the line bundles $\cm_{n-i}^V$. Under these identifications, the operator $\cd$ is a linear isomorphism between the spaces of  functions
\begin{eqnarray}\label{E:cosine-tr-compact}
\cd_i\colon C^\infty(Gr_{n-i}^V)\tilde\to C^\infty(Gr_i^V),
\end{eqnarray}
where we put subscript $i$ in $\cd_i$ for convenience. Clearly $\cd_i$ commutes with $GL(\Lam)$.

To describe $\cd_i$ more explicitly let us introduce a notation.
Let $M,N\subset V$ be linear subspaces of complementary dimension. Let us define
\begin{eqnarray}\label{E:def-s-symmetric}
s(M,N):=\mu^V((M\cap \Lam)+(N\cap \Lam)).
\end{eqnarray}
Clearly $$0\leq s(M,N)\leq 1.$$

It is easy to see that
\begin{eqnarray}\label{E:sine-another}
s(M,N)=\mu_M(p_M(N\cap \Lam))=\mu_N(p_N(M\cap\Lam)),
\end{eqnarray}
where $p_M\colon V\to V/M$, $p_N\colon V\to V/N$ are the canonical quotient maps.
Also obviously for any $T\in GL(\Lam)$ one has
\begin{eqnarray}
s(TM,TN)=s(M,N).
\end{eqnarray}

The next claim easily follows by unwinding the definitions.
\begin{claim}\label{Cl:formula-cosine}
The operator $\cd_i$ from (\ref{E:cosine-tr-compact}) is equal to
$$(\cd_i f)(M)=\int_{Gr_{n-i}^V}s(M,N) f(N)d\mu_{n-i,Haar}(N),$$
where $\mu_{n-i,Haar}$ is the $GL(\Lam)$-invariant Haar probability measure on $Gr_{n-i}^V$.
\end{claim}

\item Consider the Hermitian product of functions on Grassmannians using the Haar measures on them.
Using the symmetry in the definition of $s(M,N)$ we easily have the following adjointness property
\begin{claim}\label{Cl:adjointness}
For any $f\in  C^\infty(Gr_{n-i}^V), g\in  C^\infty(Gr_{i}^V)$ one has
$$(\cd_i f,g)=(f,\cd_{n-i}g).$$
\end{claim}
{\bf Proof.} This immediately follows from the definition (\ref{E:def-s-symmetric}) of $s(M,N)$ which is clearly symmetric with respect to $M$ and $N$. \qed

\end{paragraphlist}

\section{Hard Lefschetz theorem.}\label{S:hard-Lefschetz}
\begin{paragraphlist}
\item Let $V$ be a finite dimensional vector space over a non-Archimedean local field $\FF$. Let us fix a lattice $\Lam\subset V$.
Let $GL(\Lam)$ denote the subgroup of $GL(V)$ consisting of such transformations $T$ such that $T(\Lam)=\Lam$. Thus $GL(\Lam)\simeq GL_n(\co)$.

Let us denote by $V_1\in Val_1^\infty (V)$ the  $GL(\Lam)$-invariant non-zero valuation from Section \ref{S:poincare}. The main result of this section is
\begin{theorem}\label{T:hard-lefschtz-t}
Let $0\leq i<n/2$. The operator of multiplication by $V_1^{n-2i}$
$$Val_i^\infty(V)\to Val_{n-i}^\infty(V)$$
is an isomorphism of vector spaces.
\end{theorem}

\item The proof of the last theorem occupies the rest of this section. Since operator $\cd$ is an isomorphism, every
$\phi\in Val^\infty_i(V)$ has the form
\begin{eqnarray}\label{E:val-phi-representation}
\phi(E)=\int_{Gr_{n-i}^V}p_{E,F}^*(f(F)) \,\,\mbox{ for all }E\in Gr_i(V),
\end{eqnarray}
where $f\in C^\infty(Gr_{n-i}^V,\cm_{n-i}^V)$, and $p_{E,F}\colon E\to V/F$ is the natural map.

Let $k\leq n-i$. By Proposition \ref{P:powers-mean-width} $V_1^k\ne 0$ and is clearly $GL(\Lam)$-invariant. Hence it can be presented as
$$V_1^k(E)=\int_{Gr_{n-k}^V}p_{E,L}^*(g(L)),$$
where $g\in C^\infty(Gr_{n-k}^V,\cm_{n-k}^V)$ is $GL(\Lam)$-invariant.

It will be helpful to make some identifications induced by the lattice $\Lam\subset V$.
The choice of $\Lam$ induces a probability Haar measure $\mu_{p,Haar}$ on each Grassmannian $Gr_p^V$ and hence a trivialization of the line bundle $|\ome_p|$ of densities.
It also induces trivialization of the line bundles $\cm_p',\cl_q^V,$ and $\cm_p^V=\cm_p'\otimes |\ome_p|$ as follows.
 By the definition, the fiber of $\cm_p'$ over a subspace $H\in Gr_p^V$ is the space of Lebesgue measures $D(V/H)$.
 Let us choose the (only) Lebesgue measure $\mu_H$ on  $V/H$ such that its value on the lattice $\Lam/(\Lam\cap H)\subset V/H$ is equal to 1. Then the sections $f,g$ get identified with continuous functions
on appropriate Grassmannians which will be denoted by $\hat f,\hat g$.

Similarly, we trivialize the line bundle $\cl_q^V$ over $Gr_q^V$ by choosing in each fiber $\cl_q|_G$ the Lebesgue measure $\mu^G$ on $G$ which is equal to 1 on the lattice $\Lam\cap G\subset G$.

For a linear subspace $E\subset V$, let us denote by $vol_E$ the Lebesgue measure on $E$ such that $vol_E(E\cap \Lam)=1$.

\begin{lemma}\label{L:product-by-intrin-volume}
Assume that a valuation $\phi$ is given by (\ref{E:val-phi-representation}). Then
one has $$\phi\cdot V_1^k=c_{n,i,k} (\cd_{i+k}\circ R_{n-i-k,n-i})( \hat f),$$
where $c_{n,i,k}>0$ is a constant.
\end{lemma}
{\bf Proof.} By (\ref{E:star01}) we have for any $E\in Gr_{i+k}^V$
\begin{eqnarray}\label{E:product-int-geom}
(\phi\cdot V_1^k)(E)=\left(\int_{F\in Gr_{n-i}^V}\int_{L\in Gr_{n-k}^V}\left((p_F\times p_L)^*(f\boxtimes g)\right)(E)\right) vol_E,
\end{eqnarray}
where $p_F\colon V\to V/F$ and $p_L\colon V\to V/L$ are the natural maps, and here we denote by $p_F\times p_L$ the obvious map $V\to V/F\times V/L$
(this is different from our previous convention according to which it would be the map $V\times V\to V/F\times V/L$).

Let us rewrite the integrand in the right hand side of (\ref{E:product-int-geom}) in terms of $\hat f,\hat g$. We can factorize uniquely the map
$p_F\times p_L\colon V\to V/F\times V/L$
as $$V\overset{p_{F\cap L}}{\to} V/(F\cap L)\overset{\overline{p_F\times p_L}}{\to }V/F\times V/L.$$

The pairs $(F,L)\in Gr_{n-i}^V\times Gr_{n-k}^V$ which are transversal to each other form an open subset of full measure. Hence below we will consider
only pairs of transversal subspaces.  In this case the map $$\overline{p_F\times p_L}\colon V/(F\cap L)\to V/F\times V/L$$ is an isomorphism.

With the above identifications we have
\begin{eqnarray*}
\left((p_F\times p_L)^*(f\boxtimes g)\right)(E)=\\
\hat f(F)\cdot  (\mu_F\boxtimes \mu_L)((p_F\times p_L)(\Lam\cap E))\cdot vol_E=\\
\hat f(F)\cdot\frac{(\mu_F\boxtimes \mu_L)((p_F\times p_L)(\Lam\cap E))}{(\mu_F\boxtimes \mu_L)((p_F\times p_L)(\Lam))} \left[(\mu_F\boxtimes \mu_L)((p_F\times p_L)(\Lam))\right]\cdot vol_E=\\
\hat f(F)\cdot\frac{\mu_{F\cap L}(p_{F\cap L}(\Lam\cap E))}{\mu_{F\cap L}(p_{F\cap L}(\Lam))} \left[(\mu_F\boxtimes \mu_L)((p_F\times p_L)(\Lam))\right]\cdot vol_E=\\
\hat f(F)\cdot\mu_{F\cap L}(p_{F\cap L}(\Lam\cap E)) \left[(\mu_F\boxtimes \mu_L)((p_F\times p_L)(\Lam))\right]\cdot vol_E,
\end{eqnarray*}
where the third equality follows from the fact that linear isomorphisms ($\overline{p_F\times p_L}$ in this case)
preserve the ratio of Lebesgue measures of sets; the last equality follows from the definition of $\mu_{F\cap L}$.
To abbreviate, let us denote
\begin{eqnarray}\label{E:ccc001}
c(F,L):=(\mu_F\boxtimes \mu_L)((p_F\times p_L)(\Lam)),\\\label{E:sss001}
s(F\cap L,E):=\mu_{F\cap L}(p_{F\cap L}(\Lam\cap E)).
\end{eqnarray}
Clearly
$$0\leq c(F,L), \, s(F\cap L,E)\leq 1$$
and $c(F,L)>0$ if $F$ and $L$ are transversal. Also obviously $s(F\cap L,E)$ depends only on $F\cap L$ and not on $F$ and $L$ separately.

Let us consider the compact metrizable topological space
$$\cw:=\{(F,L,W)\in Gr_{n-i}^V\times Gr_{n-k}^V\times Gr_{n-i-k}^V|\,\, W\subset F\cap L\}.$$

\begin{claim}
$\cw$ has a structure of an analytic manifold (see Section \ref{S:analytic}) such that the natural projections $\cw\to Gr_{n-i}^V, Gr_{n-k}^V, Gr_{n-i-k}^V$ are analytic, and the natural action
$GL(V)\times\cw\to \cw$ is an analytic map.
\end{claim}
{\bf Proof.} Consider the partial flag manifolds
\begin{eqnarray*}
Z_1=\{(F,W)\in Gr_{n-i}^V\times Gr_{n-i-k}^V| \, W\subset F\},\\
Z_2:=\{(L,W)\in Gr_{n-k}^V\times Gr_{n-i-k}^V| \, W\subset L\}.
\end{eqnarray*}
$Z_1,Z_2$ are analytic manifolds by Theorem \ref{T:serre}. We have the natural analytic map
$$T\colon Z_1\times Z_2\to (Gr_{n-i-k}^V)^2$$
given by $\left((F,W_1), (L,W_2)\right)\overset{T}{\mapsto} (W_1,W_2)$. Its differential is onto at every point. Let $\Delta\subset (Gr_{n-i-k}^V)^2$ be the diagonal submanifold.
The implicit function theorem (see Section \ref{S:analytic}, paragraph \ref{item:implicit-function}) implies that $\cw=T^{-1}(\Delta)$ is an analytic submanifold of $Z_1\times Z_2$. The claimed properties of $\cw$ follow easily. \qed

\hfill

Let us continue proving Lemma \ref{L:product-by-intrin-volume}.
Clearly if $(F,L,W)\in \cw$ and if $F,L$ are transversal then $W=F\cap L$. Hence
the natural map $\cw\to Gr_{n-i}^V\times Gr_{n-k}^V$ forgetting $W$ is an isomorphism between open subsets of full measure where $E$ and $L$ are transversal.

Let $m$ be the push-forward under the inverse map of the product of the Haar probability measures
(with respect to the natural action of $GL(\Lam)$) on $Gr_{n-i}^V$ and $Gr_{n-k}^V$ with subsequent extension by 0 to the whole $\cw$.
Clearly $m$ is a $GL(\Lam)$-invariant probability measure whose value on any open subset of $\cw$ is positive (the action of $GL(\Lam)$ is diagonal). Then we can rewrite
\begin{eqnarray*}
(\phi\cdot V_1^k)(E)=\int_{(F,L,W)\in \cw}\hat f(F) c(F,L)s(W,E)dm
\end{eqnarray*}

Let us consider the analytic manifold
$$\cx:=\{(F,W)\in Gr_{n-i}^V\times Gr_{n-i-k}^V|\,\, W\subset F\}$$
and the natural $GL(\Lam)$-equivariant map $\tau\colon \cw\to \cx$  given by $\tau(F,L,W)=(F,W)$, where $GL(\Lam)$ acts diagonally.
Let $\tilde m:=\tau_*(c\cdot m)$ denote the push-forward of the measure $c\cdot m$. Clearly $\tilde m$ is $GL(\Lam)$-invariant and positive  (Haar)
measure on the $GL(\Lam)$-homogeneous space $\cx$. Then
\begin{eqnarray*}
(\phi\cdot V_1^k)(E)=\left(\int_{(F,W)\in \cx}\hat f(F) s(W,E)d\tilde m\right) vol_E=\\
c_{n,i,k}\left(\int_{W\in Gr_{n-i-k}^V}dW s(W,E)\int_{F\supset W}dF \hat f(F)\right) vol_E ,
\end{eqnarray*}
where $c_{n,i,k}>0$ is a constant, the inner integral is taken with respect to the set of all subspaces $F$ containing $W$ which can be identified with the Grassmannian $Gr_k^{V/W}$.
 The measures $dW$ on Grassmannian $Gr_{n-i-k}^V$ is the $GL(\Lam)$-invariant probability Haar measure, and the measure $dF$ is the Haar probability measure $Gr_k^{V/W}$ invariant 
 under the stabilizer of $W$ in $GL(\Lam)$.

The inner integral $\int_{F\supset W}dF \hat f(F) $ is the Radon transform $(R_{n-i-k,n-i}\hat f)(W)$. Thus we get in this notation
\begin{eqnarray}\label{E:Radon-cosine}
(\phi\cdot V_1^k)(E)=c_{n,i,k}\left(\int_{W\in Gr_{n-i-k}^V} (R_{n-i-k,n-i}\hat f)(W) \cdot s(W,E)dW\right)vol_E.
\end{eqnarray}
By Claim \ref{Cl:formula-cosine} the last expression is equal to $c_{n,i,k}(\cd_{i+k}(R_{n-i-k,n-i}\hat f))(E)$ where $c_{n,i,k}>0$.
Lemma \ref{L:product-by-intrin-volume} is proven.

\item Let us finish the proof of Theorem \ref{T:hard-lefschtz-t}. Let us assume now that $i<n/2$. Let $k=n-2i$.
Choice of lattice $\Lam\subset V$ induced a trivialization of the line bundle $\cl_i$ over $Gr_i^V$.
Indeed recall that the fiber of $\cl_i$ over $E\in Gr_i^V$ is $\cl_i|_E=D(E)$. But $D(E)=\CC\cdot vol_E$. Hence we will identify in this paragraph the space $Val_i^\infty(V)$ with a
subspace of locally constant $\CC$-valued functions on $Gr_i^V$.

   Let $\phi\in Val^\infty_i(V)$. By Lemma \ref{L:product-by-intrin-volume} we have
$$\phi\cdot V_1^{n-2i}=c_{n,i}(\cd_{n-i}\circ R_{i,n-i})(\hat f),\,\,\, c_{n,i}>0.$$


Since $R_{i,n-i}$ and $\cd_{n-i}$ are isomorphisms, the result follows. \qed



\end{paragraphlist}


\section{Fourier transform commutes with exterior product.}\label{S:fourier-ext-product}
\begin{paragraphlist}
\item In this section we fix a non-Archimedean local field $\FF$. All vector spaces $X,Y,Z,...$ will be over this field. The main result of this section is
\begin{theorem}\label{T:fourier-ext-prod-commut}
Let $X,Y$ be finite dimensional vector spaces. Let $\phi\in Val(X),\,\psi\in Val^\infty(Y)$. Then
$$\FF(\phi\boxtimes \psi)=\FF(\phi)\boxtimes \FF(\psi).$$
\end{theorem}
\begin{remark}
For convex valuations an analogue of this result was conjectured by the author in \cite{alesker-barcelona} and proved recently by Faifman and Wannerer \cite{faifman-wannerer}.
\end{remark}
We will need two lemmas to prove Theorem \ref{T:fourier-ext-prod-commut}.

\item Let $X$ and $Y=X\oplus Z$ be finite dimensional vector spaces over field $\FF$.
Fix a positive Lebesgue measure $vol_Z\in D(Z)$. Let $vol_Z^{-1}\in D(Z)^*=D(Z^\vee)$ be the corresponding Lebesgue measure
on $Z^\vee$ as defined in Lemma \ref{L:dual-measure-non-arch}.
\begin{lemma}\label{L:push-for-lemma-imb}
Let $F\colon X\inj X\oplus Z$ be the obvious imbedding given by $F(x)=(x,0)$. Let  $\phi_0\in Val(X),\, \mu\in D(X)^*$. Then
$$F_*(\phi_0\otimes \mu)=(\phi_0\boxtimes vol_Z)\otimes(\mu\otimes vol_Z^{-1}).$$
\end{lemma}
{\bf Proof.} We may assume that $\mu=vol_X^{-1}$, where $vol_X$ is a non-vanishing Lebesgue measure on $X$. By (\ref{E:yyyyy}) we have
\begin{eqnarray*}
\phi_0\boxtimes vol_Z=(F_*\otimes Id_{D(X\oplus Z)})((\phi_0\otimes vol_X^{-1})\otimes(vol_X\otimes vol_Z))=\\
F_*\phi\otimes (vol_X\otimes vol_Z).
\end{eqnarray*}
This is equivalent to the required equality. \qed

\item \begin{lemma}\label{L:fourier-upp-back-density}
Let a vector space $X$ be a direct sum $X=Z\oplus Z_1$. Let $p\colon X\to Z$ be the obvious projection. Let $\mu_Z\in D(Z)$. Then
$$\FF(p^*\mu_Z)=(\chi_{Z^\vee}\boxtimes vol_{Z_1}^{-1})\otimes (\mu_Z\otimes vol_{Z_1}),$$
where $vol_{Z_1}$ is an arbitrary non-vanishing Lebesgue measure on $Z_1$, and $vol_{Z_1}^{-1}\in D(Z_1^\vee)$ is the inverse Lebesgue measure on $Z_1^\vee$ defined in Lemma \ref{L:dual-measure-non-arch}.
\end{lemma}
{\bf Proof.} We have
\begin{eqnarray*}
\FF(p^*\mu_Z)=p^\vee_*(\FF\mu_Z)\overset{(\ref{E:fourier-density})}{=}\\
p^\vee_*(\chi_{Z^\vee}\otimes\mu_Z)\overset{\mbox{Lemma }\ref{L:push-for-lemma-imb}}{=}\\
(\chi_{Z^\vee}\boxtimes vol_{Z_1}^{-1})\otimes (\mu_Z\otimes vol_{Z_1}).
\end{eqnarray*}
\qed

\item Let us prove a special case of Theorem \ref{T:fourier-ext-prod-commut}.
\begin{lemma}\label{L:special-case-fourier-ext-prod}
Let $\mu_X\in D(X)$, $\psi\in Val(Y)$. Then
$$\FF(\mu_X\boxtimes\psi)=\FF(\mu_X)\boxtimes \FF(\psi).$$
\end{lemma}
{\bf Proof.} Since both sides are linear and continuous with respect to $\psi\in Val(Y)$, we may and will assume that
$$\psi=p^*\mu_M,$$
where $Y=M\oplus L$, $p\colon Y\to M$ is the obvious projection, $\mu_M\in D(M)$. Then we have
\begin{eqnarray*}
\FF(\mu_X\boxtimes p^*\mu_M)\overset{\mbox{Lemma }\ref{L:ext-dens-pull-back-dens}}{=}\FF((Id_X\times p)^*(\mu_X\boxtimes \mu_M))=\\
(Id_{X^\vee}\times p^\vee)_*(\FF(\mu_X\boxtimes \mu_M))\overset{(\ref{E:fourier-density})}{=}\\
(Id_{X^\vee}\times p^\vee)_*(\chi_{X^\vee\oplus M^\vee}\otimes (\mu_X\boxtimes\mu_M)).
\end{eqnarray*}
 By Lemma \ref{L:push-for-lemma-imb} we can continue
\begin{eqnarray}\label{E:8881}
\FF(\mu_X\boxtimes p^*\mu_M)=(\chi_{X^\vee\oplus M^\vee}\otimes vol_L^{-1})\otimes((\mu_X\boxtimes\mu_M)\otimes vol_L).
\end{eqnarray}
On the other hand we have
\begin{eqnarray*}
\FF(\mu_X)\boxtimes \FF(p^*\mu_M)\overset{(\ref{E:fourier-density})}{=}\\
(\chi_{X^\vee}\otimes \mu_X)\boxtimes \FF(p^*\mu_M)\overset{Lemma \ref{L:fourier-upp-back-density}}{=}\\
(\chi_{X^\vee}\otimes \mu_X)\boxtimes\left((\chi_{M^\vee}\boxtimes vol_L^{-1})\otimes(\mu_M\otimes vol_L)\right)=\\
(\chi_{X^\vee}\boxtimes \chi_{M^\vee}\boxtimes vol_L^{-1})\otimes(\mu_X\otimes \mu_M\otimes vol_L)\overset{\mbox{Prop. } \ref{P:ext-prod-1}}{=}\\
(\chi_{X^\vee\oplus M^\vee}\otimes vol_L^{-1})\otimes((\mu_X\boxtimes\mu_M)\otimes vol_L)\overset{(\ref{E:8881})}{=}\\
\FF(\mu_X\boxtimes p^*\mu_M).
\end{eqnarray*}
\qed

\item {\bf Proof of Theorem \ref{T:fourier-ext-prod-commut}.}
It suffices to assume that $\phi=p^*_M\mu$ where $p_M\colon X\to X/M, \, \mu\in D(X/M)$.
Since $\psi$ is smooth, it can be presented in the form
$\psi=\int_{N\in Gr^Y_{n-j}}p_N^*\nu(N),$
where $p_N\colon Y\to Y/N$ is the quotient map. Then we have
\begin{align}
\FF(\phi\boxtimes\psi)=\FF(p_M^*\mu\boxtimes \int_{N\in Gr^Y_{n-j}}p_N^*\nu(N))\overset{\mbox{Lemma }\ref{P:pull-back-exterior-product}}{=}\notag\\
\FF\left((p_M\times Id_Y)^*(\mu\boxtimes \int_{Gr^Y_{n-j}} p_N^*\nu(N))\right)=\notag\\
 \int_{N\in Gr^Y_{n-j}} \FF\left((p_M\times Id_Y)^*(\mu\boxtimes  p_N^*\nu(N))\right)\overset{\mbox{Lemma }\ref{L:ext-dens-pull-back-dens}}{=}\notag\\
 \int_{N\in Gr^Y_{n-j}} \FF\left(\left((p_M\times Id_Y)^*\circ (Id_{X/M}\times p_N)^*\right)(\mu\boxtimes \nu(N))\right)\overset{\mbox{Thm }\ref{T:pull-back-thm}(2)}{=}\notag\\
 \int_{N\in Gr^Y_{n-j}} \FF\left((p_M\times p_N)^*(\mu\boxtimes \nu(N))\right).\notag
\end{align}


We will rewrite the last expression with the use of Lemma \ref{L:Four-onto-comp}.
Let us fix lattices $\Lam_X\subset X$ and $\Lam_Y\subset Y$.
This choice induces positive Lebesgue measures $vol_X$ and $vol_Y$ such that $vol_X(\Lam_X)=vol_Y(\Lam_Y)=1$ and similarly on all their subspaces and quotient spaces and their dual spaces. Thus, for a subspace $M\subset X$
the corresponding Lebesgue measure $vol_M$ satisfies $vol_M(\Lam_M)=1$ where $\Lam_M=\Lam_X\cap M$, and the Lebesgue measure on $vol_{X/M}$ satisfies $vol_{X/M}(\Lam_{X/M})=1$ where $\Lam_{X/M}\subset X/M$
is the image of $\Lam_X$ under the canonical quotient map $X\to X/M$. Lebesgue measure $vol_{X^\vee}$ on the dual space $X^\vee$ is defined similarly using the dual lattice $\Lam_X^\vee\subset X^\vee$.

Then $\mu=\tilde\mu \cdot vol_{X/M}$, where $\tilde\mu\in \CC$. The choice of $\Lam_Y$ also induces an identification of the section $\nu(N)$ with $\tilde\nu(N) vol_{Y/N}$, where $\tilde\nu(N)\in \CC$. Then Lemma \ref{L:Four-onto-comp} implies that
\begin{eqnarray}\label{E:Four-ppp}
\FF(\phi\boxtimes\psi)=\tilde\mu \int_{N\in Gr^Y_{n-j}}dN \tilde\nu(N)\cdot  (q_M\times q_N)^*(vol_{M^\vee}\boxtimes vol_{N^\vee})\otimes(vol_X\otimes vol_Y),
\end{eqnarray}
where the notation is as follows. The maps $q_M\colon X^\vee\to M^\vee,\, q_N\colon Y^\vee\to N^\vee$ are the canonical quotient maps.
$vol_{X/M}\in D(X/M)$ and $vol_{Y/N}\in D(Y/N)$. Then
$$\mu=\tilde\mu\cdot vol_{X/M} \mbox{ with } \tilde\mu\in \CC.$$
Next, $dN$ is a Haar measure on $Gr^Y_{n-j}$ also induced by the lattice $\Lam_Y\subset Y$.

Let us now compute $\FF\phi\boxtimes \FF\psi$. We have
\begin{eqnarray*}
\FF\phi\boxtimes\FF\psi=\FF(p_M^*\mu)\boxtimes \FF\left(\int_{N\in Gr^Y_{n-j}}p^*_N(\nu(N))\right)\overset{Lemma \ref{L:Four-onto-comp}}{=}\\
\left[q_M^* (vol_{M^\vee})\otimes (vol_{M^\vee}^{-1}\otimes\mu)\right]\boxtimes \left[\int_{N\in Gr^Y_{n-j}}\FF(p^*_N(\nu(N)))\right]\overset{Lemma \ref{L:Four-onto-comp}}{=}\\
\left[q_M^* (vol_{M^\vee})\otimes (\underset{=\tilde\mu vol_X}{\underbrace{vol_{M^\vee}^{-1}\otimes\mu}})\right]\boxtimes\left[\int_{N\in Gr^Y_{n-j}}q^*_N(vol_{N^\vee})\otimes(\underset{\tilde \nu(N) vol_Y dN}{\underbrace{vol_{N^\vee}^{-1}\otimes\nu(N)}})\right]=\\
\left[q_M^* (vol_{M^\vee})\otimes (\tilde\mu vol_X)\right]\boxtimes \left[\int_{N\in Gr^Y_{n-j}}dN\tilde\nu(N)\cdot q_N^*(vol_{N^\vee})\otimes vol_Y\right]=\\
\left[q_M^* (vol_{M^\vee})\boxtimes \left(\int_{N\in Gr^Y_{n-j}}dN\tilde\nu(N)\cdot q_N^*(vol_{N^\vee})\right)\right]\otimes (\tilde\mu\cdot vol_X\otimes vol_Y)\overset{Prop. \ref{P:pull-back-exterior-product}}{=}\\
(q_M\times Id)^*\left[vol_{M^\vee}\boxtimes \left(\int_{N\in Gr^Y_{n-j}}dN\tilde\nu(N)\cdot q_N^*(vol_{N^\vee})\right)\right]\otimes (\tilde\mu\cdot vol_X\otimes vol_Y)=\\
(q_M\times Id)^*\left[\int_{N\in Gr^Y_{n-j}}dN\tilde\nu(N)\left(vol_{M^\vee}\boxtimes  q_N^*(vol_{N^\vee})\right)\right]\otimes (\tilde\mu\cdot vol_X\otimes vol_Y)=\\
\left[\int_{N\in Gr^Y_{n-j}}dN\tilde\nu(N)(q_M\times Id)^*\left(vol_{M^\vee}\boxtimes  q_N^*(vol_{N^\vee})\right)\right]\otimes (\tilde\mu\cdot vol_X\otimes vol_Y)\overset{Prop. \ref{P:pull-back-exterior-product}}{=}\\
\left[\int_{N\in Gr^Y_{n-j}}dN\tilde\nu(N)(q_M\times Id)^*(Id\times q_N)^*\left(vol_{M^\vee}\boxtimes  vol_{N^\vee}\right)\right]\otimes (\tilde\mu\cdot vol_X\otimes vol_Y)\overset{Thm\, \ref{T:pull-back-thm}(2)}{=}\\
\left[\int_{N\in Gr^Y_{n-j}}dN\tilde\nu(N)(q_M\times q_N)^*\left(vol_{M^\vee}\boxtimes  vol_{N^\vee}\right)\right]\otimes (\tilde\mu\cdot vol_X\otimes vol_Y).
\end{eqnarray*}
The last expression is equal to (\ref{E:Four-ppp}). \qed


\end{paragraphlist}


\section{Convolution on smooth valuations.}\label{S:convolution}
\begin{paragraphlist}
\item The goal of this section is to define a convolution on valuations and prove its basic properties. First note that if
$$\phi\in Val^\infty(X)\otimes D(X)^*,\,\, \psi\in Val^\infty(Y)\otimes D(Y)^*$$
then $\phi\boxtimes\psi$ is well defined as an element of $Val(X\times Y)\otimes D(X\times Y)^*=Val(X\times Y)\otimes D(X)^*\otimes D(Y)^*$.

Let $a\colon V\times V\to V$ be the addition map, i.e. $a(x,y)=x+y$. Let us define the convolution
$$\ast\colon (Val^\infty(V)\otimes D(V)^*)\times (Val^\infty(V)\otimes D(V)^*)\to Val^\infty(V)\otimes D(V)^*$$
by
\begin{eqnarray}\label{E:def-convolut}
\phi\ast\psi:=a_*(\phi\boxtimes\psi).
\end{eqnarray}

\begin{proposition}\label{P:convolution-smooth}
Convolution of smooth valuations is smooth
\end{proposition}
{\bf Proof.} Indeed the convolution can be considered as a linear map
$$\ast\colon (Val^\infty(V)\otimes D(V)^*)\otimes (Val^\infty(V)\otimes D(V)^*)\to Val^\infty(V)\otimes D(V)^*.$$
Evidently it commutes with the natural action of the group $GL(V)$. Hence smooth vectors are mapped into smooth ones. \qed

\item \begin{proposition}\label{P:convol-vs-product}
Let $\phi,\psi\in Val^\infty(V)$. Then
$$\FF\phi\ast\FF\psi=\FF(\phi\cdot\psi).$$
\end{proposition}
{\bf Proof.} Let $a \colon V^\vee\times V^\vee \to V^\vee$ denote the addition map. Then the dual map
$$a^\vee\colon V\to V\times V$$
is the diagonal map, i.e. $a^\vee(v)=(v,v)$. Then by definition of the convolution, push-forward, and the product we have
\begin{eqnarray*}
\FF\phi\ast\FF\psi=a_*(\FF\phi\boxtimes \FF\psi)\overset{\mbox{Thm }\ref{T:fourier-ext-prod-commut}}{=}
a_*(\FF(\phi\boxtimes\psi))=\FF\left((a^\vee)^*(\phi\boxtimes\psi)\right)=\\
\FF(\phi\cdot \psi).
\end{eqnarray*}
\qed

\item Let $vol_V\in D(V)$ be a non-zero Lebesgue measure on $V$. We denote by $vol_V^{-1}\in D(V)^*$ the element of the dual space of $D(V)$ whose
value on $vol_V$ is equal to 1. Note that $vol_V\otimes vol_V^{-1}\in D(V)\otimes D(V)^*$ is independent of a choice of $vol_V$.

Denote $n:=\dim V$.
\begin{theorem}\label{T:convolution-properties}
1) $Val^\infty(V)\otimes D(V)^*$ equipped with the convolution is a commutative associative algebra with the unit element $vol_V\otimes vol_V^{-1}$.
\newline
2) $$(Val_{n-i}^\infty(V)\otimes D(V)^*)\ast(Val_{n-j}^\infty(V)\otimes D(V)^*)\subset Val_{n-i-j}^\infty(V)\otimes D(V)^*.$$
\newline
3) The Poincar\'e duality is satisfied: the bilinear map
$$Val_i^\infty(V)\otimes D(V)^*\times Val_{n-i}^\infty(V)\otimes D(V)^*\to Val_0(V)\otimes D(V)^*=D(V)^*$$
given by $(\phi,\psi)\mapsto \phi\ast\psi$
is a perfect pairing, i.e. for any non-zero $\phi\in Val_i^\infty(V)\otimes D(V)^*$ there is $\psi\in Val_{n-i}^\infty(V)\otimes D(V)^*$
such that $\phi\ast\psi\ne 0$.
\newline
4) The hard Lefschetz type result is true: Let us fix a lattice $\Lam\subset V$. This induces an isomorphism $D(V)^*\simeq \CC$. Let
$V_{n-1}\in Val_{n-1}(V)$ be the only (up to a proportionality) $GL(\Lam)$-invariant non-zero element. Let $n/2<i\leq n$. Then the linear map
$$Val_i^\infty(V)\otimes D(V)^*\to Val_{n-i}^\infty(V)\otimes D(V)^*$$
given by $\phi\mapsto \phi\ast \underset{2i-n}{\underbrace{V_{n-1}\ast\dots\ast V_{n-1}}}$
is an isomorphism.
\end{theorem}

{\bf Proof} immediately follows from Proposition \ref{P:convol-vs-product}, the corresponding properties of the product, and the obvious fact that $\FF(V_1)$ is  proportional
to $V_{n-1}$ (after the appropriate identifications induced by the choice of lattice $\Lam\subset V$ are applied). \qed

\end{paragraphlist}

\section{Valuations invariant under a subgroup.}\label{S:subgroup}
\begin{paragraphlist}
\item Convex valuations invariant under various subgroups of $GL_n(\RR)$ have a number of interesting properties and found applications in integral geometry, see e.g. \cite{alesker-jdg-03},
\cite{bernig-aig}, \cite{bernig-faifman}, \cite{bernig-fu-annales}, \cite{fu-aig}. For a compact subgroup $G\subset GL_n(\RR)$ the author showed \cite{alesker-2000}
that the space of $G$-invariant convex valuations is finite dimensional if and only if $G$ acts transitively on the unit sphere, and in this case all $G$-invariant valuations are smooth \cite{alesker-04}.
Furthermore in this case  the algebra $Val^G(\RR^n)$ of $G$-invariant valuations equipped either with product or convolution satisfies Poincar\'e duality, hard Lefschetz theorem, and Hodge-Riemann bilinear relations inherited from $Val^\infty(\RR^n)$.

It is well known that any compact subgroup $G\subset GL_n(\FF)$ is conjugated to a subgroup of $GL_n(\co)$. In particular any maximal compact subgroup of $GL_n(\FF)$ is conjugated to $GL_n(\co)$.

We will see below that in the non-Archimedean case there are a lot of compact subgroups $G\subset GL(V)$ such that the subspace of $G$-invariant valuations is finite dimensional and all such valuations are smooth.

\item Let $V$ be an $n$-dimensional vector space over a non-Archimedean local field $\FF$.
The group $GL(V)$ has many subgroups which are simultaneously open and compact, more precisely such subgroups form a basis of neighborhoods of $I_n\in GL(V)$.
\begin{proposition}\label{P:G-inv}
Let $G\subset GL(V)$ be an open and compact subgroup. Then the space $Val^G(V)$ of $G$-invariant valuations is finite dimensional and all its elements are smooth, i.e.
$$Val^G(V)\subset Val^\infty(V).$$
\end{proposition}
{\bf Proof.} Orbits of any open subgroup of $GL(V)$ on $Gr_i^V$ are open. Indeed the map $GL(V)\to Gr_i^V$ given by $g\mapsto g(E_0)$ is a submersion for any $E_0\in Gr_i^V$.
In particular $G$-orbits are open. Since $Gr_i^V$ is compact and different orbits are disjoint, there are finitely many of them. By definition of a smooth valuation
$$Val^G(V)\subset Val^\infty(V).$$
\qed

\item Obviously
$$Val^G(V)=\oplus_{i=0}^n Val_i^G(V).$$
It is easy to see that $(Val^G(V),\cdot)$ is a finite dimensional graded subalgebra of $(Val^\infty(V),\cdot)$ satisfying Poincar\'e duality and hard Lefschetz theorem (with respect to $V_1\in Val_1^G(V)$ which is invariant under a maximal compact subgroup containing $G$).

Similarly $(Val^G(V)\otimes D(V)^*,\ast)$ is a finite dimensional graded subalgebra of $(Val^\infty(V)\otimes D(V)^*,\ast)$ satisfying Poincar\'e duality and hard Lefschetz (with respect to $V_{n-1}$ which is invariant under a maximal compact subgroup containing $G$).
(Note that the compatibility with the grading is given by Theorem \ref{T:convolution-properties}(2).)

The Fourier transform establishes an isomorphism of algebras
$$\FF\colon (Val^G(V),\cdot)\tilde\to (Val^G(V^\vee)\otimes D(V),\ast).$$

\end{paragraphlist}

\end{document}